\documentclass[aop,MSNbibl,seceqn,citesort,dvips]{arximspdf}

\doi{10.1214/10-AOP636}
\volume{40}
\issue{2}
\pubyear{2012}
\firstpage{648}
\lastpage{694}

\makeatletter

\newtheorem{theorem}{Theorem}[section]
\newtheorem{lemma}[theorem]{Lemma}
\newtheorem{cor}[theorem]{Corollary}
\newtheorem{prop}[theorem]{Proposition}

\newproclaim{rmk}[theorem]{Remark}
\newproclaim{defn}[theorem]{Definition}
\newproclaim{standing}[theorem]{Standing Assumptions}

\newcommand{\ccdots}{\cdots}

\newcommand{\al}{\alpha}
\newcommand{\oooverline}{\hspace*{1pt}\overline}

\newcommand{\cA}{\mathcal{A}}
\newcommand{\cB}{\mathcal{B}}
\newcommand{\cE}{\mathcal{E}}
\newcommand{\F}{\mathcal{F}}

\newcommand{\I}{\mathcal{I}}
\newcommand{\LL}{\mathcal{L}}

\newcommand{\cS}{\mathcal{S}}
\newcommand{\cZ}{\mathcal{Z}}

\newcommand{\E}{\mathbb{E}_{\alpha}}
\newcommand{\PPP}{\mathbb{P}_{\alpha}}

\newcommand{\R}{\mathbb{R}}
\newcommand{\N}{\mathbb{N}}
\newcommand{\Z}{\mathbb{Z}}
\newcommand{\din}{d_{\mathrm{in}}}
\newcommand{\dout}{d_{\mathrm{out}}}
\newcommand{\dmax}{d_{\mathrm{max}}}

\newcommand{\h}{\mathbf{h}}

\newcommand{\per}{\operatorname{per}}
\newcommand{\tr}{\operatorname{tr}}

\newcommand{\Per}{\operatorname{Per}}
\newcommand{\Sp}{\operatorname{Sp}}
\newcommand{\Var}{\operatorname{Var}}

\makeatother

\begin{document}
\begin{frontmatter}

\title{Random subshifts of finite type}
\runtitle{Random subshifts of finite type}

\begin{aug}
\author[A]{\fnms{Kevin} \snm{McGoff}\corref{}\ead[label=e1]{mcgoff@math.umd.edu}\ead[label=u1,url]{www.math.umd.edu/\textasciitilde mcgoff}}
\runauthor{K. McGoff}
\affiliation{University of Maryland}
\address[A]{Department of Mathematics \\
University of Maryland \\
College Park, Maryland 20742-4015\\
USA\\
\printead{e1}\\
\printead{u1}} 
\end{aug}

\received{\smonth{6} \syear{2010}}
\revised{\smonth{10} \syear{2010}}

%
\begin{abstract}
Let $X$ be an irreducible shift of finite type (SFT) of positive
entropy, and let $B_n(X)$ be its set of words of length $n$. Define
a~random subset $\omega$ of $B_n(X)$ by independently choosing each word
from~$B_n(X)$ with some probability $\al$. Let $X_{\omega}$ be the
(random) SFT built from the set $\omega$. For each $0 \leq\al\leq1$
and $n$ tending to infinity, we compute the limit of the likelihood
that $X_{\omega}$ is empty, as well as the limiting distribution of
entropy for $X_{\omega}$. For $\al$ near~$1$ and $n$ tending to
infinity, we show that the likelihood that $X_{\omega}$ contains a
unique irreducible component of positive entropy converges
exponentially to~$1$. These results are obtained by studying certain
sequences of random directed graphs. This version of ``random SFT''
differs significantly from a previous notion by the same name, which
has appeared in the context of random dynamical systems and bundled
dynamical systems.
\end{abstract}

%
\begin{keyword}[class=AMS]
\kwd[Primary ]{37B10}
\kwd[; secondary ]{37B40}
\kwd{60C05}.
\end{keyword}
\begin{keyword}
\kwd{Random subshifts of finite type}
\kwd{entropy}.
\end{keyword}

\end{frontmatter}


\section{Introduction}\label{intro}

A shift of finite type (SFT) is a dynamical system defined by finitely
many local transition rules. These systems have been studied for their
own sake \cite{Kitchens,LM}, and they have also served as important
tools for understanding other dynamical systems \cite{Keane,Bowen,DGS}.

Each SFT can be described as the set of bi-infinite sequences on a
finite alphabet that avoid a finite list of words over the alphabet.
Thus, there are only countably many SFTs up to the naming of letters in
an alphabet.

For the sake of simplicity, we state our results in terms of SFTs in
the \hyperref[intro]{Introduction}, even though we prove more general results in terms of
sequences of directed graphs in the subsequent sections. Let $X$ be a
nonempty SFT (for definitions, see Section \ref{SFTpresentations}). Let
$B_n(X)$ be the set of words of length $n$ that appear in $X$. For $\al
$ in $[0,1]$, let $\PPP$ be the probability measure on the power set of
$B_n(X)$ given by choosing each word in $B_n(X)$ independently with
probability $\al$. The case $\al= 1/2$ puts uniform measure on the
subsets of $B_n(X)$. For notation, let $\Omega_n$ be the power set of
$B_n(X)$. To each subset $\omega$ of $B_n(X)$, we associate the SFT
$X_{\omega}$ consisting of all points $x$ in $X$ such that each word of
length $n$ in $x$ is contained in $\omega$. With this association, we
view $\PPP$ as a probability measure on the SFTs $X_{\omega}$ that can
be built out of the subsets of $B_n(X)$. Briefly, if $X$ has entropy $\h
(X) = \log\lambda>0$ and $n$ is large, then a typical random SFT
$X_{\omega}$ is built from about $\al\lambda^n$ words, an $\al$
fraction of all the words in $B_n(X)$, but not all of these words will
occur in any point in $X_{\omega}$.

Our main results can be stated as follows. Let $\zeta_X(t)$ denote the
Artin--Mazur zeta function of $X$ (see Definition \ref{zetaDefn}). The
first theorem deals with the likelihood that a randomly chosen SFT is empty.
\begin{theorem} \label{EmptyThmIntro}
Let $X$ be a nonempty SFT with entropy $\h(X) = \log\lambda$. Let
$\mathcal{E}_n \subset\Omega_n$ be the event that $X_{\omega}$ is
empty. Then for $\al$ in $[0,1]$,
\[
\lim_{n \rightarrow\infty} \PPP(\mathcal{E}_{n}) = \cases{(\zeta_X(\al
))^{-1}, &\quad if $\al\in[0, 1/\lambda)$,\cr
0, &\quad if $\al\in[1/\lambda,1]$.}
\]
\end{theorem}

Thus, when $\al$ is in $[0,1/\lambda)$, there is an asymptotically
positive probability of emptiness. The next theorem gives more
information about what happens when $\al$ lies in $[0, 1/\lambda)$.
\begin{theorem} \label{SubCritThmIntro}
Let $X$ be a nonempty SFT with entropy $\h(X) = \log\lambda$. Let
$\mathcal{Z}_n \subset\Omega_n$ be the event that $X_{\omega}$ has
zero entropy, and let $I_n$ be the random variable on $\Omega_n$ which
is the number of irreducible components of $X_{\omega}$. Then for
$0\leq\al< 1/\lambda$,
\begin{longlist}[(3)]
\item[(1)] $\lim_{n \rightarrow\infty} \PPP(\mathcal{Z}_{n}) = 1$;
\item[(2)] the sequence $(I_{n})$ converges in distribution to the random
variable $I_{\infty}$ such that $\mathbb{P}(I_{\infty} = 0) = (\zeta
_X(\al))^{-1}$ and for $k \geq1$,
\[
\mathbb{P}(I_{\infty} = k) = (\zeta_X(\al))^{-1} \mathop{\sum_{S
\subset\N}}_{|S|=k} \prod_{s \in S} \frac{ \al^{|\gamma_{s}|}}{1-\al
^{|\gamma_{s}|}},
\]
where $\{\gamma_i \}_{i=1}^{\infty}$ is an enumeration of the periodic
orbits in $X$;
\item[(3)] the random variable $I_{\infty}$ has exponentially decreasing
tail and therefore finite moments of all orders.
\end{longlist}
\end{theorem}

Our next result describes the entropy of the typical random SFT when
$\al$ lies in $(1/\lambda, 1]$.
\begin{theorem} \label{EntropyThmIntro}
Let $X$ be an SFT with positive entropy $\h(X) = \log\lambda$. Then
for $1/ \lambda< \al\leq1$ and $\varepsilon>0$,
\[
\lim_{n \rightarrow\infty} \PPP\bigl( |\h(X_{\omega}) - \log(\al\lambda)|
\geq\varepsilon\bigr) = 0,
\]
and the convergence to this limit is exponential in $n$.
\end{theorem}

Finally, we have a result concerning the likelihood that a random SFT
will have a unique irreducible component of positive entropy when $\al$
is near~$1$.

\begin{theorem} \label{MSFTthmIntro}
Let $X$ be an irreducible SFT with positive entropy\break $\h(X) = \log
\lambda$. Let $W_{n} \subset\Omega_n$ be the event that $X_{\omega}$
has a unique irreducible component~$C$ of positive entropy and~$C$ has
the same period as $X$. Then there exists $c>0$ such that for $1-c <
\al\leq1$,
\[
\lim_{n \rightarrow\infty} \PPP( W_{n} ) = 1;
\]
furthermore, the convergence to this limit is exponential in $n$.
\end{theorem}

There have been studies of other objects called random subshifts of
finite type in the literature \cite
{BogGund,BogGund2,FS,Kifer2,Kifer3,Kifer4,Kifer,KL}, but the objects
studied here are rather different in nature. The present work is more
closely related to perturbations of SFTs, which have already appeared
in works by Lind \cite{Lind2} in dimension~$1$ and by Pavlov \cite
{Pavlov} in higher dimensions. In those works, the main results
establish good uniform bounds for the entropy of an SFT obtained by
removing any single word of length $n$ from a sufficiently mixing SFT
as $n$ tends to infinity. Random SFTs may also be interpreted as
dynamical systems with holes
\cite{CM1,CM2,CMT1,CMT2,CvdB,Demers1,Demers2,DWY,DY,LMD,LopesM}, in
which case the words of length~$n$ in $X$ that are forbidden in the
random SFT $X_{\omega}$ are viewed as (random) holes in the original
system $X$. The question of whether an SFT defined by a~set of
forbidden words is empty has been studied in formal language theory and
automata theory, and in that context it amounts to asking whether the
set of forbidden words is \textit{unavoidable}
\cite{Bell,CHP,HS}. Also, the random SFTs considered here can be viewed
as specific instances of random matrices (see~\cite{BS,Mehta}) or
random graphs (see \cite{ABS,BNP,Dur,ER1,ER2,Gr1,Gr2,NP}), and the
concept of directed percolation on finite graphs has appeared in the
physics literature in the context of directed networks \cite{OP,HOR}.
To the best of our knowledge, the specific considerations that arise
for our random SFTs seem not to have appeared in any of this wider literature.

The paper is organized as follows. Section \ref{Preliminaries} contains
the necessary background and notation, as well as some preliminary
lemmas. The reader familiar with SFTs and directed graphs may prefer to
skip Sections \ref{SFTpresentations} and~\ref{SeqsGraphs}, referring
back as necessary. In Section \ref{Emptiness} we discuss the likelihood
that a~random SFT is empty, and, in particular, we prove Theorem \ref
{EmptyThmIntro}. The remainder of the main results are split into two
sections according to two cases: $\al\in[0, 1/\lambda)$ and $\al\in
(1/\lambda,1]$. The case $\al\in[0, 1/\lambda)$ is treated in Section
\ref{subcriticalPhase}, and the case $\al\in(1/\lambda,1]$ is
addressed in Section \ref{supercriticalPhase}. Section \ref{Remarks}
discusses some corollaries of the main results.

\section{Preliminaries} \label{Preliminaries}

\subsection{Shifts of finite type and their presentations} \label
{SFTpresentations}

For a detailed treatment of SFTs and their presentations, see\vadjust{\goodbreak} \cite
{LM}. In this section we describe three ways to present an SFT: with a
finite list of forbidden words over a finite alphabet, with a finite,
directed graph, or with a square, nonnegative integer matrix.

Let $\mathcal{A}$ be a finite set, which we will call the
\textit{alphabet}. An element $b \in\cA^{n}$ is called a \textit{word} of
length $n$. Let $\Sigma= \mathcal{A}^{\mathbb{Z}}$, endowed with the
product topology induced by the discrete topology on $\mathcal{A}$.
Then $\Sigma$ is a compact metrizable space, which is called the
\textit{full shift} on $\cA$. Let $\sigma\dvtx\Sigma\rightarrow\Sigma$ be the
left shift, that is, for $x = ( x_{i} )$ in~$\Sigma$, let
$(\sigma(x))_{i} = x_{i+1}$. With this definition $\sigma$ is a~homeomorphism of $\Sigma$.

A subset $X$ of $\Sigma$ is called shift-invariant if $\sigma(X) = X$.
A closed, shift-invariant subset of $\Sigma$ is called a
\textit{subshift}. For any subshift $X$, the \textit{language}~$\mathcal
{B}(X)$ of $X$ is the collection of all finite words (blocks) that
appear in some sequence $x$ in $X$. Note that $\cB(X) = \bigcup B_n(X)$,
where $B_n(X)$ is the set of all words of length $n$ that appear in
some sequence $x$ in $X$. [By convention we set $B_0(X) = \{\varepsilon
\}$, where $\varepsilon$ denotes the empty word.] Given a set $\mathcal
{F}$ of words on $\cA$, we may define a subshift $X(\mathcal{F})$ as
the set of sequences $x$ in~$\Sigma$ such that no word in $\mathcal{F}$
appears in $x$. One may check that this procedure indeed defines a
subshift. If $X$ is a subshift and there exists a \textit{finite} set
of words $\mathcal{F} = \{F_1,\ldots, F_k \}$ such that $X = X(\mathcal
{F})$, then $X$ is called a \textit{subshift of finite type} (SFT).

The natural notion of isomorphism for SFTs is called conjugacy. Two
SFTs $X$ and $Y$ are \textit{conjugate}, written $X \cong Y$, if there
exists a homeomorphism $\phi\dvtx X \to Y$ such that $\phi\circ\sigma=
\sigma\circ\phi$. An SFT $X$ is \textit{irreducible} if for every two
nonempty open sets $U$ and $V$ and every $N$ in $\N$, there exists $n
\geq N$ such that $\sigma^n(U) \cap V \neq\varnothing$. An SFT $X$ is
\textit{mixing} if for every two nonempty open sets $U$ and $V$ in $X$,
there exists $n_0$ in $\N$ such that for all $n \geq n_0$, we have
$\sigma^n(U) \cap V \neq\varnothing$. Mixing and irreducibility are
conjugacy-invariant. We now define the higher block presentations of an SFT.
\begin{defn} \label{higherBlockDef}
Let $X$ be an SFT. The $n$-\textit{block presentation} of $X$, denoted
$X^{[n]}$, is defined as follows. The alphabet for $X^{[n]}$ is
$B_n(X)$. We define the code $\phi_n\dvtx X \to B_n(X)^{\Z}$ by the equation
%
\begin{equation}
\phi_n(x)_i = x[i,i+n-1]
\end{equation}
for all $x$ in $X$. Then $X^{[n]} = \phi_n(X)$. For all $n \geq1$, we
have that $X^{[n]} \cong X$, where the conjugacy is given by $\phi_n$.
\end{defn}
\begin{defn}
The \textit{entropy} of an SFT $X$ is defined as $\h(X) = \lim_n \frac
{1}{n}\times {\log}|B_n(X)|$.
\end{defn}

Alternatively, one may define SFTs in terms of finite directed graphs.
A~directed graph $G = (V,E)$ consists of a set of vertices $V$ and a
set of~ed\-ges~$E$ such that for each edge $e \in E$, there is a unique
initial vertex, \mbox{$i(e) \in V$}, and a unique terminal vertex, $t(e) \in
V$. We view the edge $e$ as going from~$i(e)$ to~$t(e)$.\vadjust{\goodbreak} We allow
self-loops, but for the sake of convenience we assume (without loss of
generality for our considerations) that there are no multiple edges. In
this paper we make the standing convention that ``graph'' means
directed graph. We will collect our standing assumptions in Standing
Assumptions~\ref{StAssump}.
\begin{defn} \label{EdgeShiftDefn}
Given a directed graph $G$, we define the \textit{edge shift} $X_G$ to
be the set of all bi-infinite (oriented) walks on $G$,
\textit{that is}, $X_G = \{ x \in E^{\Z}\dvtx t(x_j) = i(x_{j+1}) \mbox{ for all } j
\in\Z\}$.
\end{defn}

Any edge shift is an SFT (trivially). Let us show that any SFT is
conjugate to an edge shift. If $X = X(\F)$ is an SFT and $\F$ is a
finite set of forbidden words, then $X \cong X_G$, where $G = (V,E)$ is
defined as follows. Let $n_0 = \max\{ |F|\dvtx F \in\F\}$. Then let $V
= B_{n_0-1}(X)$ and $E = B_{n_0}(X)$. Further, for any edge $e \in
B_{n_0}(X)$, we let $i(e) = e[1,n_0-1]$ and $t(e) = e[2,n_0]$. The same
construction works with $n$ in place of $n_0$ for any $n \geq n_0$.

If $G$ is a graph such that $X \cong X_G$, we say that $X_G$ is an edge
presentation of $X$, or sometimes just a presentation of $X$. The
\textit{adjacency matrix} $A$ of a~directed graph $G$ may be defined as
follows. Fix an enumeration of the vertices in $G$. Then let $A_{k \ell
}$ be the number of distinct edges $e$ in $G$ such that $i(e) = v_k$
and $t(e) = v_{\ell}$. A square, nonnegative integral matrix $A$ is
\textit{irreducible} if for each pair $i,j$ and each $N$, there exists
$n >N$ such that $(A^n)_{ij} >0$. A~matrix $A$ is \textit
{nondegenerate} if it has no zero row and no zero column. If $A$ is
nondegenerate, then the edge shift $X_{G}$ is irreducible if and only
if $A$ is irreducible. Also, if $A$ is nondegenerate, then the edge
shift~$X_G$ is mixing if and only if there exists $n_0$ such that for
all $n \geq n_0$ and all pairs $i,j$, it holds that $(A^n)_{ij} > 0$. A
matrix is \textit{primitive} if it satisfies the latter property. A~%
path in $G$ is a finite sequence $\{e_j\}_{j=1}^n$ of edges such that
$t(e_{j}) = i(e_{j+1})$ for $j = 1,\ldots, n-1$. If $b = b_1\ccdots b_n$
is a path in $G$, we say that $b$ goes from vertex~$i(b_1)$ to vertex
$t(b_n)$. We denote by $B_k(G)$ the set of paths of length $k$ in $G$.
By convention, we set $B_0(G) = V$.
\begin{defn} For a path $b$ in $G$, let $V(b)$ and $E(b)$ be the set of
vertices and the set of edges traversed by $b$, respectively.
\end{defn}
\begin{defn}
Let $X$ be an SFT. An \textit{irreducible component} $Y$ of $X$ is a
nonempty, maximal SFT contained in $X$ such that $Y$ is irreducible.
Let~$G$ be a graph. An \textit{irreducible component} $C$ of $G$ is a
nonempty, maximal subgraph of $G$ such that the adjacency matrix of $C$
is irreducible. The reader should be advised that in some papers the
definition of irreducible component includes trivial components (a
single vertex with no edges adjacent to it), but the definition given
here does not include trivial components.
\end{defn}
\begin{defn} \label{nBlockGraphDefn}
Let $G$ be a finite, directed graph. For $n \geq1$, define $G^{[n]} =
(V^{[n]},E^{[n]})$, the $n$-\textit{block graph} of $G$, as
follows.\vspace*{1pt}
Let $V^{[n]} = B_{n-1}(G)$\vadjust{\goodbreak} and $E^{[n]} = B_{n}(G)$, such that if $e
\in E^{[n]}$, then $i(e) = e[1,n-1]$ and $t(e) = e[2,n]$. Note
that\vspace*{1pt}
$G^{[1]} = G$.
\end{defn}

If $X = X_{G}$ for some graph $G$, then it follows immediately from the
definitions that $X^{[n]} = X_{G^{[n]}}$.
\begin{defn} \label{pthPowerDefn}
Let $G = (V,E)$ be a graph. For $p$ in $\N$, we define the $p$th power
graph, $G^p = (V^p,E^p)$, as follows. Let $V^p = V$ and $E^p = B_p(G)$.
If $b=b_1 \ccdots b_p$ is an edge in $G^p$, then we let $i(b) = i(b_1)$
and $t(b) = t(b_1)$.
\end{defn}
\begin{defn} \label{TransposeDefn}
Let $G = (V,E)$ be a graph. Define the transpose graph, $G^T =
(V^T,E^T)$, as follows. Let $V^T = V$ and $E^T = E$, where an edge $e$
in $G^T$ goes from $t(e)$ to $i(e)$. In other words, the transpose
graph is just the graph formed by reversing the direction of all the
edges in $G$.
\end{defn}

Given a square, nonnegative, integral matrix $A$, one may also define
an SFT $X_A$ as follows. Let $G$ be a directed graph whose adjacency
matrix is exactly $A$ (such a graph always exists). Then let $X_A$ be
the edge shift defined by $G$.

Recall\vspace*{1pt} the following basic facts (which may be found in \cite{LM}). For
an SFT $X$, we have $\h(X) = {\inf_n \frac{1}{n} \log}|B_n(X)|$. If $X$
is a nonempty SFT and $X = X_A$ for a square, nonnegative integral
matrix $A$, then $\h(X) = \log\lambda$, where $\lambda$ is the
spectral radius of $A$. By the Perron--Frobenius theorem, if $A$ is
nonnegative and irreducible, then there exists a strictly positive
(column) vector $v$ such that $Av = \lambda v$, and there exists a
strictly positive (row) vector~$w$ such that $w A = \lambda w$.
Furthermore, $v$ and $w$ are each unique up to a positive scalar.
\begin{defn} \label{nonzeroSpec} For any nonnegative integer matrix
$A$, let $\lambda_A$ be the spectral radius of $A$, and let $\chi_A$ be
the characteristic polynomial of $A$. Then let $\Sp_{\times}(A)$ be the
\textit{nonzero spectrum} of the matrix $A$, which is defined as the
multiset of nonzero roots of $\chi_A$ listed according to their
multiplicity. If~$A$ is the adjacency matrix of the graph $G$, we
define $\lambda_G = \lambda_A$ and $\Sp_{\times}(G) = \Sp_{\times}(A)$.
\end{defn}

If $X_{A} \cong X_{B}$ for two nonnegative integral matrices $A$ and
$B$, then $\Sp_{\times}(A) = \Sp_{\times}(B)$. Also, if $A$ is
primitive, then $\max\{ |\beta|\dvtx \beta\in\Sp_{\times}(A) \setminus
\{\lambda_A\} \} < \lambda_A$. Finally, if $A$ is irreducible, then
there exists a unique $\sigma$-invariant Borel probability measure $\mu
$ on $X_{A}$ of maximal entropy. Let us describe some basic properties
of $\mu$. We associate a word $b = b_1 \ccdots b_k$ in $X$ to the
cylinder set $C_{b} = \{ x \in X\dvtx x[1,k] = b\}$. In this way we
interpret the measure of words in~$\cB(X)$ as the measure of the
corresponding cylinder set. Let $v$ be a positive right eigenvector of
$A$ and $w$ a positive left eigenvector of $A$, and suppose they are
normalized so that $w \cdot v = 1$. Our standing assumption that there
are no multiple edges means that $A_{ij} \leq1$ for all $i,j$. Then
for a vertex $u$ in~$V$, we have $\mu(u)=w_u v_u$, and for $b \in
B_n(X_A)$, we have that
%
\begin{equation} \label{measureOfAPath}
\mu(b) = w_{i(b_1)} \lambda_A^{-n} v_{t(b_n)}.
\end{equation}

Now we define two objects, the period and the zeta function, which
contain combinatorial information about the cycles in a graph $G$
(alternatively, one may refer to the periodic points in an SFT $X$).
\begin{defn} \label{periodDefn} For an SFT $X$, let $\per(X)$ be the
greatest common divisor of the sizes of all periodic orbits in $X$. For
a graph $G$, let $\per(G)$ be the greatest common divisor of the
lengths of all cycles in $G$.
\end{defn}
\begin{defn} \label{zetaDefn}
Let $X$ be an SFT and $N_p = |\{x \in X\dvtx \sigma^p(x) = x\}|$.
Then the Artin--Mazur zeta function of $X$ (see \cite{LM}) is, by definition,
\[
\zeta_{X}(t) = \exp\Biggl( \sum_{p=1}^{\infty} \frac{N_p}{p} t^p \Biggr).
\]
For a graph $G$, let $\zeta_G = \zeta_{X_G}$.
\end{defn}

For a graph $G$, note that $|\{x \in X_G\dvtx \sigma^p(x) = x\}|$ is the
number of cycles of (not necessarily least) period $p$ in $G$, and
\[
\zeta_{G}(t) = \frac{1}{\det(I-tA)} = \prod_{\lambda\in\Sp_{\times
}(G)} \frac{1}{1-\lambda t}.
\]
Also, $\zeta_G$ has radius of convergence $1/\lambda_G$ and $\lim_{t
\to1/\lambda_G^{-}} \zeta_G(t) = +\infty$.

\subsection{Sequences of graphs under consideration} \label{SeqsGraphs}

In this work we consider sequences of graphs $(G_n)$ that grow in some
way. A particular example of such a sequence is the sequence of
$n$-block graphs of an SFT $X$. Indeed, by taking $(G_n)$ to be such a
sequence in Theorems \ref{EmptyThm}, \ref{SubCritThm}, \ref{EntropyThm}
and \ref{MSFTthm}, we obtain the theorems stated in the \hyperref[intro]{Introduction}.
Generalizing to the graph setting also allows one to consider sequences
of graphs presenting SFTs which are conjugate to a fixed SFT $X$, where
the sequences need not be the $n$-block sequence for $X$. To indicate
the generality of the arguments further, though, we formulate and prove
the results for sequences of graphs that do not necessarily present
conjugate SFTs. Before we move on to these results, we need to define
several notions regarding the manner of growth of the sequence $(G_n)$.

Let $G$ be a finite, directed graph with adjacency matrix $A$. We will
have use for the following notation.
\begin{defn} \label{PerDefn}
Let
\[
\Per_p(G) = \{ b \in B_p(G)\dvtx i(b_1) = t(b_p)\}\quad \mbox{and}\quad \Per(G) =
\bigcup_{p=1}^{\infty} \Per_p(G).
\]
For $b$ in $\Per_p(G)$, let $\theta(b)$ be the set of all paths $c$ in
$\Per_p(G)$ such that there exists a natural number $\ell$ such that $c
= b_{\tau^{\ell}(1)} \ccdots b_{\tau^{\ell}(p)}$, where $\tau$ is the
permutation of $\{1,\ldots, k\}$ defined in cycle notation by $(1 \ccdots k)$.
\end{defn}
\begin{defn}
For each vertex $u$ in $G$, let $\dout(u) = |\{ e \in E\dvtx i(e) = u\}|$
and $\din(u) = |\{ e \in E\dvtx t(e) = u\}|$. Then let
\[
\dmax(G) = \max\{ \max(\dout(u), \din(u))\dvtx u \in V\}.
\]
\end{defn}

In order to measure the separation of periodic orbits in $G$, we make
the following definition.
\begin{defn} \label{Defnofz} Let
\[
z(G) = \max\Biggl\{ n \geq0\dvtx \forall b, c \in\bigcup_{p=1}^{n} \Per_p(G)
\mbox{ with } c \notin\theta(b), V(b) \cap V(c) = \varnothing\Biggr\},
\]
where $V(b)$ is the set of vertices traversed by the path $b$.
\end{defn}

As a measure of the size of $G$, we consider the following quantity.
\begin{defn} If $A$ has spectral radius $\lambda>1$, then let
\[
m(G) = \lceil{\log_{\lambda} }|V| \rceil.
\]
\end{defn}

To measure a range for uniqueness of paths in $G$, we make the
following definitions.
\begin{defn} \label{UDefn} Let
\begin{eqnarray*}
U_1(G) &=& \sup\{ n\dvtx \forall i,j \mbox{ it holds that } (A^n)_{ij}
\leq1 \}, \\
U_2(G) &=& \sup\bigl\{ n\dvtx \forall u \in V \mbox{ and } 1 \leq s < t \leq
n,\\
&&\hphantom{\sup\bigl\{ } |\{b \in B_t(X)\dvtx i(b_1)=u, b_s = b_t\}| \leq1 \bigr\}, \\
U(G) &=& \min(U_1(G),U_2(G)).
\end{eqnarray*}
\end{defn}

We use the transition length as a type of diameter of $G$.
\begin{defn} Let
\[
R(G) = \inf\{ n\dvtx \forall i,j, \exists k \leq n, (A^k)_{ij} >0 \}.
\]
\end{defn}

Here we briefly recall the notion of the weighted Cheeger constant of
an irreducible, directed graph $G$. The weighted Cheeger constant was
defined\vadjust{\goodbreak} and studied in \cite{Chung}. Let $\mu$ be the measure of
maximal entropy of $X_G$, and let $F\dvtx E \rightarrow[0,1]$ be given by
$F(e) = \mu(e)$. For any vertex $v$ in $V$, let $F(v) = \sum_{i(e) = v}
F(e) = \sum_{t(e) = v} F(e)$. Then for any subset of vertices $S
\subseteq V$, let $F(S) = \sum_{v \in S} F(v)$, and for any two subsets
$S,T \subseteq V$, let
\[
F(S,T) = \mathop{\sum_{i(e)\in S}}_{t(e) \in T} F(e).
\]
In general, $F(S,T)$ is not symmetric in $S$ and $T$ since $G$ is
directed. Let~$E(S,\allowbreak T)$ be the set of edges $e$ in $G$ such that $i(e)
\in S$ and $t(e) \in T$. Let $\overline{S} = V \setminus S$.

\begin{defn} \label{CheegerDefn} The weighted Cheeger constant of $G$
is defined as
\[
c_w(G) = \inf_{\varnothing\subsetneq S \subsetneq V} \frac
{F(S,\overline{S})}{\min(F(S),F(\overline{S}))},
\]
and the unweighted Cheeger constant of $G$ is defined as
\[
c(G) = \inf_{0 < |S| \leq|V|/2} \frac{|E(S,\overline{S})|}{|S|}.
\]
\end{defn}
\begin{defn} \label{uniformExpSeqDef}
We say that $G$ is a directed $b$-expander graph if $c(G) \geq b$.
Also, a sequence of directed graphs $(G_n)$ is a \textit{uniform
expander sequence}, if there exists a $b>0$ such that $G_n$ is a
directed $b$-expander for each $n$.
\end{defn}

We will also have use for the following quantity related to the
spectral gap of~$G$.
\begin{defn} \label{SpGapDefn}
Let $g(G) = \min\{ 1- \frac{|\lambda_i|}{\lambda}\dvtx \lambda_i \in\Sp
_{\times}(G) \setminus\{\lambda\} \}$.
\end{defn}

We make the following standing assumptions, even though some of the
statements we make may hold when these restrictions are relaxed. In
particular, Theorems \ref{EmptyThm} and \ref{SubCritThm} do not require
that $A_n$ is irreducible, nor do they require that $\lambda>1$ (see
Remark \ref{relaxRmk}).
\begin{standing} \label{StAssump} Recall that ``graph'' means
directed\break
graph. Let $(G_n)$ be a sequence of graphs with an associated sequence
of adjacency matrices $(A_n)$. Unless otherwise stated, we will make
the following assumptions:
\begin{itemize}
\item for each $n$, each entry of $A_n$ is contained in $\{0,1\}$;
\item each $A_n$ is irreducible;
\item for each $n$, $\Sp_{\times}(A_n) = \Sp_{\times}(A_1)$;
\item$\lambda:= \lambda_{A_1} >1$;
\item$\lim_n m(G_n) = \infty$.\vadjust{\goodbreak}
\end{itemize}
\end{standing}
\begin{rmk}
Note that $|{\Per_p}(G_n)| = \tr(A_n^p)$, which depends only on $\Sp
_{\times}(A_n)$ and $p$. Therefore, the standing assumptions imply that
$|{\Per_p}(G_n)|$ does not depend on $n$, and, therefore, $\per(G_n)$ and
$ \zeta_{G_n}$ do not depend on $n$.
\end{rmk}

Additional conditions that we place on sequences of graphs will come
from the following list. [Different theorems will require different
assumptions, but the sequence of $n$-block graphs of an irreducible
graph with spectral radius greater than $1$ will satisfy conditions
(C1)--(C8) below by Proposition \ref{higherGraphSeqLemma}.]

\begin{defn} \label{conditions} We define the following conditions on a
sequence of graphs $(G_n)$ with a sequence of adjacency matrices $(A_n)$:
\begin{longlist}[(C8)]
\item[(C1)] there exists $\Delta>0$ such that $\dmax(G_n) \leq\Delta$ for
all $n$ (bounded degree);
\item[(C2)] $z(G_n)$ tends to infinity as $n$ tends to infinity (separation
of periodic points);
\item[(C3)] there exists $C > 0$ such that $z(G_n) \geq C m(G_n)$ for all $n$
(fast separation of periodic points);
\item[(C4)] there exists $C > 0$ such that $U(G_n) \geq m(G_n)- C$ for all
$n$ (local uniqueness of paths);
\item[(C5)] there exists $C >0$ such that $R(G_n) \leq m(G_n) +C$ for all $n$
(small diameter);
\item[(C6)] there exists $K >0$ such that $\max_{u \in V_n} \mu(u) \leq K \min
_{u \in V_n} \mu(u)$ for all~$n$ (bounded distortion of vertices) and
$\max_{e \in E_n} \mu(e) \leq K \min_{e \in E_n} \mu(e)$ for all $n$
(bounded distortion of edges);
\item[(C7)] there exists $K >0$ such that $\max_i w_i^n \leq K \min_i w_i^n$
and $\max_i v_i^n \leq K \min_i v_i^n$ for all $n$, where $w^n$ is a
positive left eigenvector of $A_n$ and $v^n$ is a positive right
eigenvector of $A_n$ (bounded distortion of weights);
\item[(C8)] $(G_n)$ is a uniform expander sequence, and $(G_n^T)$ is a
uniform expander sequence (forward/backward expansion).
\end{longlist}
\end{defn}

Now we establish some lemmas, which will be used in the subsequent sections.
\begin{lemma} \label{C7impliesC1C6}
Let $(G_n)$ be a sequence of graphs satisfying the Standing Assumptions
\ref{StAssump}. Then \textup{(C7)} implies \textup{(C1)} and \textup{(C6)}
for both $(G_n)$ and~$(G_n^T)$.
\end{lemma}
\begin{pf}
First note that if (C7) holds for~$(G_n)$, then it also holds for~$(G_n^T)$ since a positive left eigenvector for $A_n^T$ is given by
$(v^n)^T$ and a positive right eigenvector for $A_n^T$ is given by
$(w^n)^T$. Therefore, we only need to show that~(C7) for $(G_n)$
implies (C1) and (C6) for $(G_n)$ [since the same argument will apply
to $(G_n^T)$].

Let $w^n$ and $v^n$ be positive left and right eigenvectors for $A_n$,
respectively, and assume\vadjust{\goodbreak} that $w^n \cdot v^n = 1$. Recall with this
normalization, if $u$ is a vertex in~$V_n$, then $\mu(u) = w^n_u
v^n_u$. Then condition (C7) implies that there exists $K>0$ such that
for all $n$,
\begin{eqnarray*}
\max_u \mu(u) &\leq& \max_u w^n_u \max_u v^n_u \leq K^2 \min_u w^n_u
\min_u v^n_u \\
&\leq& K^2 \min_u w^n_u v^n_u = K^2 \min_u \mu(u).
\end{eqnarray*}
Similarly, (C7) implies that there exists $K'>0$ such that for all $n$,
we have that $\max_{e\in E_n} \mu(e) \leq K' \min_{e\in E_n} \mu(e)$
[recall that $\mu(e) = w^n_{i(e)} \lambda^{-1} v^n_{t(e)}$]. Thus,~(C7)
implies (C6).

Note that for $e$ in $E_n$, we have that
\[
\mu(e | i(e)) = \frac{ w^n_{i(e)} \lambda^{-1} v^n_{t(e)} }{ w^n_{i(e)}
v^n_{i(e)}} = \frac{v^n_{t(e)}}{\lambda v^n_{i(e)}}.
\]
Then condition (C7) implies that there exists a uniform constant $K>0$
such that $\mu(e | i(e)) \geq K^{-1}$ for all $n$ and all $e$ in $E_n$.
We also have that
\[
\mu(u) = \sum_{e\dvtx i(e) = u} \mu(e) \geq\sum_{e\dvtx i(e) = u} K^{-1} \mu
(u) = |\{ e\dvtx i(e) = u\}| K^{-1} \mu(u).
\]
Since $G_n$ is irreducible (by Standing Assumptions \ref{StAssump}), we
know that $\mu(u) >0$, and, therefore, we have that for any $n$, and
any $u$ in $V_n$,
\[
|\{ e \in E_n\dvtx i(e) = u\}| \leq K,
\]
which implies that $\max_u \dout(u)$ is uniformly bounded in $n$. A
similar argument shows that $\max_u \din(u)$ is uniformly bounded in
$n$, which shows that $\dmax(G_n)$ is uniformly bounded in $n$ and
gives (C1).
\end{pf}

Recall that for a graph $G$, the quantities $g(G)$ and $c_w(G)$ were
defined in Definitions \ref{SpGapDefn} and \ref{CheegerDefn}, respectively.
\begin{lemma} \label{mixExpLemma}
Let $G$ be a graph with primitive adjacency matrix $A$. Then it holds
that $c_{w}(G) \geq\frac{1}{2}g$.
\end{lemma}
\begin{pf}
This lemma is a consequence of \cite{Chung}, Theorems 4.3 and 5.1, as
we now explain. Since $A$ is primitive, there exists a strictly
positive vector~$v$ and $\lambda\geq1$ such that $Av = \lambda v$.
Let $P$ be the stochastic matrix defined\vspace*{1pt} by $P_{ij} = \frac{A_{ij}
v_j}{\lambda v_i}$. Then $P$ is the transition probability matrix
corresponding to the random walk defined by the measure of maximal
entropy $\mu$ on $X_G$. We have\vspace*{1pt} that $\Sp_{\times}(P) = \frac{1}{\lambda
} \Sp{\times}(A)$. Given such a transition probability matrix, Chung
defines a Laplacian $L$ and proves~(\cite{Chung}, Theorem 4.3), that the
smallest nonzero eigenvalue of $L$, denoted~$\lambda_1$, satisfies the
following inequality:
%
\begin{equation} \label{Chung4.3}
\min\bigl\{1- |\rho|\dvtx \rho\in\Sp_{\times}(P) \setminus\{1\} \bigr\} \leq
\lambda_1.
\end{equation}
We remark that the left-hand side of the inequality in \cite{Chung},
Theorem 4.3, is equal\vadjust{\goodbreak} to the left-hand side of (\ref
{Chung4.3}) since $A$ is primitive (not just irreducible). Note that
the left-hand side of (\ref{Chung4.3}) equals $g(G)$, as
defined in Definition~\ref{SpGapDefn}. After defining the weighted
Cheeger constant (as in Definition \ref{CheegerDefn}), Chung proves
(\cite{Chung},~Theorem~5.1), that
%
\begin{equation} \label{Chung5.1}
c_w(G) \geq\tfrac{1}{2} \lambda_1.
\end{equation}
Combining the inequalities in (\ref{Chung4.3}) and (\ref
{Chung5.1}), we obtain the desired inequality.\vspace*{-1pt}
\end{pf}

Recall that the $p$th power graph was defined in Definition
\ref{pthPowerDefn}.\vspace*{-1pt}
\begin{lemma} \label{irrExpLemma}
Let $G$ be a graph with an irreducible adjacency matrix. Let $p = \per
(G)$. Let $G^{p,0}$ be an irreducible component of $G^p$, the $p$th
power graph of $G$. Let $g = g(G^{p,0})$ (which does not depend on the
choice of irreducible component in $G^p$). Then there exists $b >0$,
depending only on $g$ and $p$, such that $c_w(G) \geq b$.\vspace*{-1pt}
\end{lemma}
\begin{pf}
Let $G$, $p$ and $g$ be as in the statement of the lemma. If $p = 1$,
then Lemma \ref{mixExpLemma} immediately gives the result. Now we
assume $p \geq2$. The fact that $G$ is irreducible and $\per(G) = p$
implies that there is a partition of the vertices into $p$ nonempty
subsets, $V = \bigcup_{j=0}^{p-1} V^j$, such that for each edge $e$ with
$i(e) \in V^j$, it holds that $t(e) \in V^{j+1}$, where the
superscripts are taken modulo $p$. Let $X = X_G$ (Definition \ref
{EdgeShiftDefn}), and for each $j = 0,\ldots, p-1$, let $X_j = \{x \in
X\dvtx i(x_0) \in V^j\}$. For any set $S \subset V$ with $0 < |S| < |V|$
and $j = 0,\ldots, p-1$, define
\begin{eqnarray*}
C_S & = & \{x \in X\dvtx i(x_0) \in S\},\qquad \overline{C}_S = X_G \setminus C_S,
\\
C_S^j & = & X_j \cap C_S \quad\mbox{and}\quad \overline{C}{}^j_S = X_j \cap\overline{C}_S.
\end{eqnarray*}
Recall that we denote by $\mu$ the measure of maximal entropy on $X$,
and we may write $c_w(G)$ as follows:
\begin{eqnarray*}
c_w(G) & = & \inf_{\varnothing\subsetneq S \subsetneq V} \frac{ \mu(C_S
\cap\sigma^{-1} \overline{C}_S )}{ \min( \mu(C_S), \mu(\overline
{C}_S)) } \\
& = & \inf_{\varnothing\subsetneq S \subsetneq V} \max\biggl(
\frac{ \mu(C_S \cap\sigma^{-1} \overline{C}_S )}{ \mu(C_S) }, \frac{
\mu(C_S \cap\sigma^{-1} \overline{C}_S )}{ \mu(\overline{C}_S) }\biggr).
\end{eqnarray*}
We also use the following notation:
%
\begin{equation} \label{r_iDefn}
r_i = \frac{\mu(C_S^i)}{\mu(C_S)}\quad \mbox{and}\quad \overline{r}_i = \frac{\mu
(\overline{C}{}^i_S)}{\mu(\overline{C}_S)}.
\end{equation}

Let us establish a useful inequality. For $i = 0,\ldots, p-1$ and $1
\leq\ell\leq p$, note that each point $x$ in $C_S^i \cap\sigma^{-\ell
}\oooverline{C}{}^{i+\ell}_S$ also lies in $C_S^j \cap\sigma^{-1}\oooverline
{C}{}^{j+1}_S$ for $j = \min\{ k > 0\dvtx \sigma^k x \notin C_S\}$. Thus,
%
\begin{equation} \label{CSestimate}
\mu( C_S^i \cap\sigma^{-\ell}\oooverline{C}{}^{i+\ell}_S ) \leq\sum
_{j=0}^{p-1} \mu(C_S^j \cap\sigma^{-1} \oooverline{C}{}^{j+1}_S) = \mu(C_S
\cap\sigma^{-1} \oooverline{C}_S).\vadjust{\goodbreak}
\end{equation}

To complete the proof, we will find $b>0$ in terms of $g$ and $p$ so
that for $S \subset V$ with $0 < |S| < |V|$, we have that
%
\begin{equation} \label{SboundEqn}
b \leq\max\biggl( \frac{ \mu(C_S \cap\sigma^{-1} \oooverline{C}_S )}{ \mu
(C_S) }, \frac{ \mu(C_S \cap\sigma^{-1} \oooverline{C}_S )}{ \mu
(\overline{C}_S) }\biggr).
\end{equation}
The bound $b$ will be the minimum of four bounds, each coming from
a~particular type of set $S \subset V$.

Consider the following conditions on the set $S$, which we will use to
break our proof into cases:
\begin{longlist}[(III)]
\item[(I)] there exists $i \in\{ 0,\ldots, p-1\}$ such that $\mu(C_S^i) \in
\{0,1\}$;
\item[(II)] $\mu(C_S^i) \leq1/2p$ for each $i$, or $\mu(C_S^i) \geq1/2p$
for each $i$;
\item[(III)] $1/4p \leq\mu(C_S^i) \leq3/4p$ for each $i$.
\end{longlist}
Now we consider cases.

\textit{Case}: (I) holds, \textit{that is}, there exists $i \in\{ 0,\ldots, p-1\}$ such that $\mu(C_S^i) \in\{0,1\}$. Assume first that $\mu
(C_S^i) = 0$, which implies that $\mu(\overline{C}{}^i_S) = \mu(X_i)$.
Choose $j$ such that $\mu(C_S^j) = \max_{k} \mu(C_S^{k})$, and finally
choose $1 \leq\ell\leq p$ such that $j+\ell= i$ ($\operatorname{mod}p$). Then by
inequality (\ref{CSestimate}) and the shift-invariance of~$\mu$, we
have that
\[
\frac{ \mu(C_S \cap\sigma^{-1} \oooverline{C}_S )}{ \mu(C_S) } \geq\frac
{ \mu(C_S^j \cap\sigma^{-\ell}\oooverline{C}{}^{j+\ell}_S )}{ \mu(C_S) }
\geq \frac{ \mu(C_S^j \cap\sigma^{-\ell} X_{j+\ell}) }{ p \max_k \mu
(C_S^k) } = \frac{ \mu(C_S^j) }{ p \mu(C_S^j) } = \frac{1}{p}.
\]
Now assume $\mu(C_S^i) = 1$. Choose $j$ such that $\mu(\overline
{C}_S^j) = \max_{k} \mu(\overline{C}{}^{k}_S)$, and finally choose $1
\leq\ell\leq p$ such that $i+\ell= j$ ($\operatorname{mod} p$). Then by (\ref
{CSestimate}) and the shift-invariance of $\mu$,
\[
\frac{ \mu(C_S \cap\sigma^{-1} \oooverline{C}_S )}{ \mu(\overline{C}_S)
} \geq\frac{ \mu(C_S^i \cap\sigma^{-\ell}\oooverline{C}{}^{i+\ell}_S )}{
\mu(\overline{C}_S) } \geq \frac{ \mu(X_i \cap\sigma^{-\ell}\oooverline
{C}_S^j) }{ p \max_k \mu(\overline{C}{}^k_S) } = \frac{ \mu(\overline
{C}_S^j) }{ p \mu(\overline{C}{}^j_S) } = \frac{1}{p}.
\]
Let $b_1 = 1/p$, and note that if condition (I) holds, then the
inequality in~(\ref{SboundEqn}) holds with $b_1$ in place of $b$.

\textit{Case}: (I) does not hold, but (II) holds, \textit{that is}, $0<
\mu(C_S^i) \leq1/2p$ for all $i$, or $1 > \mu(C_S^i) \geq1/2p$ for
all $i$. Assume first that $0< \mu(C_S^i) \leq1/2p$ for all $i$. Since
$\sum_i r_i =1$ and $r_i \geq0$ for all $i$, there exists $j$ such
that $r_j \geq1/p$. Then by (\ref{CSestimate}) and the definition of
$r_i$ in (\ref{r_iDefn}),
\begin{eqnarray*}
\frac{ \mu(C_S \cap\sigma^{-1} \oooverline{C}_S )}{ \mu(C_S) } &\geq&\frac
{ \mu(C_S^j \cap\sigma^{-p} \oooverline{C}{}^j_S )}{ \mu(C_S) } = r_j \frac
{ \mu(C_S^j \cap\sigma^{-p} \oooverline{C}{}^j_S )}{ \mu(C_S^j) }\\
&\geq&
\frac{1}{p} \cdot\frac{ \mu(C_S^j \cap\sigma^{-p} \oooverline{C}{}^j_S
)}{ \mu(C_S^j) }.
\end{eqnarray*}
Let $G^{p,j}$ be the irreducible component of $G^p$ with vertex set
$V^j$. Then $G^{p,j}$ has primitive adjacency matrix, and $g =
g(G^{p,j}) >0$. Lemma \ref{mixExpLemma} gives that $c_w(G^{p,j}) \geq
\frac{1}{2} g$. Since $\mu(C_S^j) \leq1/2p$ and $\mu(C_S^j) + \mu
(\overline{C}{}^j_S) = \mu(X_j) = 1/p$, we have that $\min(\mu(C_S^j), \mu
(\overline{C}{}^j_S)) = \mu(C_S^j)$, and, thus,
\[
\frac{ \mu(C_S^j \cap\sigma^{-p} \oooverline{C}{}^j_S )}{ \mu(C_S^j) }
\geq c_w(G^{p,j}) \geq\frac{1}{2}g.
\]
Let $b_2 = g/2p$. We have shown that for $S$ such that $\mu(C_S^i) \leq
1/2p$ for each~$i$, the inequality in (\ref{SboundEqn}) holds with
$b_2$ in place of $b$. For $S$ such that $1> \mu(C_S^i) \geq1/2p$ for
each $i$, choose $j$ such that $\overline{r}_j \geq1/p$. Then an
analogous argument gives that the inequality in (\ref{SboundEqn}) holds
with $b_2$ in place of $b$.

\textit{Case}: (III) holds, \textit{that is}, $1/4p \leq\mu(C_S^i)
\leq3/4p$ for all $i$. A simple calculation yields that $r_i \geq
1/3p$ and $\overline{r}_i \geq1/3p$ for each $i$. Using (\ref
{CSestimate}), we see that for each~$j$,
%
\begin{eqnarray} \label{caseIIIestimateI}
\frac{ \mu(C_S \cap\sigma^{-1} \oooverline{C}_S )}{ \mu(C_S) } &\geq&\frac
{ \mu(C_S^j \cap\sigma^{-p} \oooverline{C}{}^j_S )}{ \mu(C_S) } = r_j \frac
{ \mu(C_S^j \cap\sigma^{-p} \oooverline{C}{}^j_S )}{ \mu(C_S^j) } \nonumber\\[-8pt]\\[-8pt]
&\geq&
\frac{1}{3p} \cdot\frac{ \mu(C_S^j \cap\sigma^{-p} \oooverline{C}{}^j_S
)}{ \mu(C_S^j) }\nonumber
\end{eqnarray}
and
%
\begin{eqnarray} \label{caseIIIestimateII}
\frac{ \mu(C_S \cap\sigma^{-1} \oooverline{C}_S )}{ \mu(\overline{C}_S)
} &\geq&\frac{ \mu(C_S^j \cap\sigma^{-p} \oooverline{C}{}^j_S )}{ \mu
(\overline{C}_S) } = \overline{r}_j \frac{ \mu(C_S^j \cap\sigma^{-p}
\overline{C}{}^j_S )}{ \mu(\overline{C}{}^j_S) }\nonumber\\[-8pt]\\[-8pt]
&\geq&\frac{1}{3p} \cdot
\frac{ \mu(C_S^j \cap\sigma^{-p} \oooverline{C}{}^j_S )}{ \mu(\overline
{C}_S^j) }.\nonumber
\end{eqnarray}
Then since $G^{p,j}$ has a primitive adjacency matrix, Lemma \ref
{mixExpLemma} and inequalities (\ref{caseIIIestimateI}) and (\ref
{caseIIIestimateII}) give that the inequality in (\ref{SboundEqn})
holds with $b_3 := g/6p$ in place of~$b$.

\textit{Case}: each of (I), (II) and (III) does not hold, \textit{that
is}, we assume that $S$ is such that $0 < \mu(C_S^i) < 1$ for each $i$,
there exists $i_1$ and $i_2$ such that $\mu(C_S^{i_1}) > 1/2p$ and $\mu
(C_S^{i_2}) < 1/2p$, and there exists $i_3$ such that either $\mu
(C_S^{i_3}) < 1/4p$ or $\mu(C_S^{i_3}) > 3/4p$. Suppose first that $\mu
(C_S^{i_3}) < 1/4p$. Choose $j$ such that $\mu(C_S^j) = \max_k \mu
(C_S^k)$, and choose $1 \leq\ell\leq p$ such that $j + \ell= i_3$
($\operatorname{mod} p$). Calculation gives that $\mu(C_S^{i_3}) < \frac{1}{2} \mu
(C_S^j)$. Then by (\ref{CSestimate}) and the shift-invariance of $\mu$,
\begin{eqnarray*}
\frac{ \mu(C_S \cap\sigma^{-1} \oooverline{C}_S )}{ \mu(C_S) } &\geq&
\frac{ \mu(C_S^{j} \cap\sigma^{-\ell}\oooverline{C}{}^{j+\ell}_S )}{ p
\mu(C_S^j) } \geq\frac{\mu(C_S^{j}) - \mu(C_S^{i_3})}{ p \mu(C_S^j)}
\\
&\geq& \frac{\mu(C_S^j)-{(1/2)}\mu(C_S^j)}{p \mu(C_S^j)} =
\frac{1}{2p}.
\end{eqnarray*}
Now assume $\mu(C_S^{i_3}) > 3/4p$. Choose $j$ such that $\mu(C_S^j) =
\max_k \mu(C_S^k)$ and choose $1 \leq\ell\leq p$ such that $j+\ell=
i_2$ ($\operatorname{mod} p$). Calculation reveals that $\mu(C_S^{i_2}) < \frac
{2}{3}\mu(C_S^j)$. Then by (\ref{CSestimate}) and the shift-invariance
of $\mu$,
\begin{eqnarray*}
\frac{ \mu(C_S \cap\sigma^{-1} \oooverline{C}_S )}{ \mu(C_S) } &\geq&\frac
{ \mu(C_S^{j} \cap\sigma^{-\ell}\oooverline{C}{}^{j+\ell}_S )}{ p \mu
(C_S^j) } \geq\frac{\mu(C_S^{j}) - \mu(C_S^{i_2})}{ p \mu(C_S^j)} \\
&\geq& \frac{\mu(C_S^j)-(2/3)\mu(C_S^j)}{p \mu(C_S^j)} = \frac{1}{3p}.
\end{eqnarray*}
Let $b_4 = 1/3p$. We have shown that for $S$ in this case, the
inequality in (\ref{SboundEqn}) holds with $b_4$ in place of $b$.

Now let $b = \min(b_1,b_2,b_3,b_4) = \min( 1/p, g/2p, g/6p, 1/3p) =
g/6p$, which depends only on $g$ and $p$. We have shown that $c_w(G)
\geq b$.
\end{pf}

Recall that the transpose graph $G^T$ of a graph $G$ was defined in
Definition~\ref{TransposeDefn}.
\begin{lemma} \label{UnifExpLemma}
Let $(G_n)$ be a sequence of graphs satisfying the Standing Assumptions
\ref{StAssump} and such that both $(G_n)$ and $(G_n^T)$ have bounded
degrees and bounded distortion of edges and vertices [conditions
\textup{(C1)} and \textup{(C6)} in Definition~\ref{conditions}]. Then
$(G_n)$ and $(G_n^T)$ are both uniform expander sequences [condition
\textup{(C8)} in Definition \ref{conditions}].
\end{lemma}
\begin{pf}
We check that conditions (C1) and (C6) for $(G_n)$ together imply that
$(G_n)$ is a uniform expander sequence, and then the same argument will
apply to $(G_n^T)$ since (C1) and (C6) also hold for $(G_n^T)$.

Recall the following notation. Let $F\dvtx E_n \rightarrow[0,1]$ be given
by $F(e) = \mu(e)$, where $\mu$ is the measure of maximal entropy on
$X_{G_n}$. Also, $c_w(G_n)$ denotes the weighted Cheeger constant of
$G_n$ (Definition \ref{CheegerDefn}). By the Standing Assumptions~\ref
{StAssump}, $\Sp_{\times}(G_n) = \Sp_{\times}(G_1)$ for each $n$.
Therefore, $\per(G_n)$ does not depend on $n$, and we let $p = \per
(G_1)$. Let $G_n^{p,0}$ be an irreducible component of the $p$th power
graph of $G_n$, and let $g_n = g(G_n^{p,0})$. Since $g_n$ only depends
on the nonzero spectrum of $G_n$, which is constant in $n$ by the
Standing Assumptions \ref{StAssump}, we have the $g_n$ is constant in
$n$. Let $g=g_1$. By Lemma \ref{irrExpLemma}, there exists $b_n >0$,
depending only $g_n$ and $\per(G_n)$, such that $c_w(G_n) \geq b_n$.
Since we have that $g_n = g$ and $\per(G_n) = p$ for all $n$, we may
choose $b := b_1$, and we obtain that $c_w(G_n) \geq b >0$ for all $n$.

Now we relate $c_w(G_n)$ to $c(G_n)$ (Definition \ref{CheegerDefn})
using properties (C1) and~(C6). For notation, let $m = m(G_n)$. Since
$(G_n)$ satisfies conditions (C1) and~(C6), there exists $K_1, K_2 >0$
such that for every $n$ and every subset $S \subset V_n$,
\[
K_1 |S| \lambda^{-m} \leq F(S) \leq K_2 |S| \lambda^{-m}
\]
and
\[
K_1 |E_n(S,\overline{S})| \lambda^{-m} \leq F(S,\overline{S}) \leq K_2
|E_n(S,\overline{S})| \lambda^{-m}.
\]
We already have that $c_w(G_n) \geq b$, which implies that for every
$S$ such that $\varnothing\subsetneq S \subsetneq V_n$,
\[
b \leq\frac{F(S,\overline{S})}{ \min(F(S),F(\overline{S})) } \leq\frac
{K_2 |E_n(S,\overline{S})| \lambda^{-m}}{ \min(F(S),F(\overline{S})) }.
\]
Now assume $0 < |S| \leq|V_n|/2$. If $\min(F(S),F(\overline{S})) =
F(S)$, then $\min(F(S)$, $F(\overline{S})) = F(S) \geq K_1 |S| \lambda
^{-m}$. If $\min(F(S)$, $ F(\overline{S})) = F(\overline{S})$, then we have\break
$\min(F(S)$, $F(\overline{S})) = F(\overline{S}) \geq K_1 |\overline{S}|
\lambda^{-m} \geq K_1 |S| \lambda^{-m}$. Combining these estimates
gives that for all $S$ such that $0< |S| \leq|V_n|/2$, we obtain that
\[
|E_n(S,\overline{S})| \geq b \frac{K_1}{K_2} |S|,
\]
which shows that $(G_n)$ is a uniform $(b \frac{K_1}{K_2})$-expander sequence.
\end{pf}
\begin{lemma} \label{implicationsLemma}
Let $(G_n)$ be a sequence of graphs satisfying the Standing Assumptions
\ref{StAssump} and bounded distortion of weights [condition
\textup{(C7)} in Definition~\ref {conditions}]. Then:
\begin{longlist}[(1)]
\item[(1)] there exists $K >0$ such that for all $n$, $k$ and $S \subset B_k(G_n)$,
\[
K^{-1} |S| \leq\lambda^{m(G_n)+k} \mu(S) \leq|S| K;
\]
\item[(2)] there exists a $K >0$ such that for all $n$, $k$, $e \in E_n$,
and $S \subset B_k(G_n)$,
\[
K^{-1} |S \cap C_e^{n,k}| \leq\lambda^{k} \mu(S | C_e^{n,k}) \leq K |S
\cap C_e^{n,k}|,
\]
where $C_e^{n,k} = \{b \in B_k(G_n)\dvtx b_1 = e\}$;
\item[(3)] there exists $K >0$ such that for all $n$, $k$, and $1 \leq s < t
\leq k$, it holds that $\mu(A_{s,t}) \leq K \lambda^{-m(G_n)}$, where
$A_{s,t} = \{b \in B_k(G_n)\dvtx b_s = b_t\}$;
\item[(4)] there exists $K >0$ such that for all $n$, $k > U(G_n)$, and $u
\in V_n$, it holds that $\mu( \Per_k(G_n) | C_u^{n,k}) \leq K \lambda
^{-U(G_n)}$, where $C_u^{n,k} = \{ b \in B_k(G_n)\dvtx i(b_1) = u\}$ and
$U(G_n)$ was defined in Definition \ref{UDefn}.
\end{longlist}
\end{lemma}
\begin{pf} For notation, let $m = m(G_n)$ and $U = U(G_n)$.

\textit{Proof of} (1).
We have that
\[
1 = \sum_{u \in V_n} \mu(u) = \sum_{u \in V_n} w^n_u v^n_u.
\]
Then condition (C7) implies that there exists $K_1>0$ such that for
each $n$ and $u$ in~$V_n$,
\[
K_1^{-1} |V_n|^{-1} \leq w^n_u v^n_u \leq K_1 |V_n|^{-1}.
\]
By the definition of $m$, there exists $K_2 >0$ such that $K_2^{-1}
|V_n|^{-1} \leq\lambda^{-m} \leq K_2 |V_n|^{-1}$. It follows that
there exists $K_3 >0$ such that for each $n$ and $u$ in~$V_n$,
\[
K_3^{-1} \lambda^{-m} \leq w^n_u v^n_u \leq K_3 \lambda^{-m}.
\]
Then (C7) implies that there exists $K_4 >0$ such that for any $n$ and
any three vertices $u$, $u_1$ and $u_2$ in $V_n$,
\[
K_4^{-1} w^n_{u_1} v^n_{u_2} \leq w^n_u v^n_u \leq K_4 w^n_{u_1} v^n_{u_2}.
\]
Finally, we conclude that there exists $K_5 >0$ such that for each $n$,
$k$ and $b$ in $B_k(G_n)$, we have that
\[
K_5^{-1} \lambda^{-(m+k)} \leq\mu(b) = w^n_{i(b)} \lambda^{-k}
v^n_{t(b)} \leq K_5 \lambda^{-(m+k)}.
\]
The statement in (1) follows.

\textit{Proof of} (2).
The statement in (2) follows from the statement in (1) and the fact
that $\mu(C_e^{n,k}) = \mu(e)$.

\textit{Proof of} (3).
Note that from (1) we have that there exists $K>0$ such that
\[
\mu(A_{s,t}) = \sum_{\gamma\in\Per_{t-s}(G_n)} \mu(\gamma) \leq K
\lambda^{-(m+t-s)} |{\Per_{t-s}}(G_n)|.
\]
Since $\Sp_{\times}(A_n)$ does not depend on $n$ by our Standing
Assumptions \ref{StAssump}, we have that $|{\Per_{t-s}}(G_n)|$ does not
depend on $n$. Clearly, $|{\Per_{t-s}}(G_n)| \lambda^{-(t-s)}$ is bounded
as $t-s$ tends to infinity. Therefore, there exists $K'$ such that
\[
\mu(A_{s,t}) \leq K' \lambda^{-m}
\]
as desired.

\textit{Proof of} (4).
By (2), we have that there exists $K_1>0$ such that for all $n$, $k
>U$, and $u$ in $V_n$,
\[
\mu(\Per_k(G_n) | C^{n,k}_u) \leq K_1 \lambda^{-k} |{\Per_k}(G_n) \cap C_u^{n,k}|.
\]
By (2), there exists $K_2 >0$ such that for all $n$, $k > U$, and $u$
in $V_n$,
\[
|B_{k-U}(G_n) \cap C_u^{n,k-U}| \leq K_2 \lambda^{k-U}.
\]
By definition of the uniqueness parameter $U$, each path in
$B_{k-U}(G_n) \cap C_u^{n,k-U}$ can be continued in at most one way to
form a path in $\Per_k(G_n) \cap C_u^{n,k}$. Therefore, with $K_3 = K_1
K_2 >0$, we have that for all $n$, $k >U$, and $u$ in $V_n$,
\[
\mu(\Per_k(G_n) | C^{n,k}_u) \leq K_1 K_2 \lambda^{-k} \lambda^{k-U} =
K_3 \lambda^{-U}.
\]
\upqed\end{pf}
\begin{prop} \label{higherGraphSeqLemma}
Let $G_1$ be a graph with irreducible adjacency matrix~$A_1$
having
entries in $\{0,1\}$ and spectral radius $\lambda> 1$. Let $G_n =
G_1^{[n]}$ for $n \geq2$. Then the\vadjust{\goodbreak} sequence $(G_n)$ satisfies the
Standing Assumptions \ref{StAssump} and conditions \textup{(C1)--(C8)}. Moreover,
\begin{longlist}
\item$\dmax(G_n)=\dmax(G_1)$ for all $n$;
\item there exists $C >0$ such that $|m(G_n) - n| \leq C$ for all $n$;
\item$z(G_n) \geq\frac{1}{2}(n-1)$ for all $n$;
\item$U(G_n) \geq n-1$ for all $n$;
\item$R(G_n) \leq n + R(G_1)$ for all $n$.
\end{longlist}
\end{prop}
\begin{pf}
One may easily check from the definitions that each $A_n$ has entries
in $\{0,1\}$, each $A_n$ is irreducible, and $\Sp_{\times}(G_n) = \Sp
_{\times}(G_1)$. We show below that $m(G_n)$ tends to infinity as $n$
tends to infinity, which gives that~$(G_n)$ satisfies the Standing
Assumptions \ref{StAssump}.

The set of in-degrees that appear in $G_n$ is constant in $n$, and so
is the set of out-degrees that appear in $G_n$. Therefore, $\dmax(G_n)
= \dmax(G_1)$, which implies condition (C1).

By definition, $m(G_n) = \lceil\log_{\lambda} |V_n| \rceil$. Since
$G_n = G_1^{[n]}$, we have that $|V_n| = |B_{n-1}(G_1)|$. By the
standard Perron--Frobenius theory, there exist cons\-tants~$K_1$ and $K_2$
such that $K_1 \lambda^n \leq|B_n(G_1)| \leq K_2 \lambda^n$. It
follows that there exists a~constant $C >0$ such that $|m(G_n) -n| \leq
C$, and, in particular, $m(G_n)$ tends to infinity.

Recall the higher-block coding map $\phi_n\dvtx X_{G_1} \to X_{G_n}$ (see
Definition \ref{higherBlockDef}). If $x$ is a point in $X_{G_1}$, then
let $V_n(x)$ be the set of vertices in $G_n$ traversed by $\phi_n(x)$.
Let us show that $z(G_n) \geq(n-1)/2$. Recall Fine and Wilf's theorem~\cite{FW},
which can be stated as follows. Let $x$ be a periodic
sequence with period $p$, and $y$ be a periodic sequence with period
$q$. If $x[i+1,i+n] = y[i+1,i+n]$ for $n \geq p+q - \gcd(p,q)$ and $i$
in $\Z$, then $x = y$. It follows from this theorem that if $x$ and $y$
lie in distinct periodic orbits of~$X_{G_1}$ and have periods less than
or equal to $(n-1)/2$, then $V_n(x) \cap V_n(y) = \varnothing$. Thus,
$z(G_n) \geq(n-1)/2$, and, in particular, $(G_n)$ satisfies conditions
(C2) and (C3).

Note that the map $\phi_n$ gives a bijection between $B_k(G_n)$ and
$B_{k+n-1}(G_1)$ for all $k \geq0$. Using this map, we check that
$U(G_n) \geq n-1$ as follows. For any two paths $b,c \in B_{n-1}(G_1)$,
there is at most one path of length $2n-2$ in $G_1$ of the form $bc$
(since every edge in such a path is specified by either~$b$ or $c$).
This fact implies that $U_1(G_n) \geq n-1$. Now if $b$ is in
$B_{n-1}(G_1)$ and $1 \leq s < t \leq n-1$ are given, then there is at
most one path $c$ in $B_{t+n-2}(G_1)$ such that $c[1,n-1]=b$ and
$c[s,s+n-2] = c[t,t+n-2]$; indeed, if $c$ is such a path, then
$c[1,n-1]$ is determined by $b$, and $c[n,t+n-1]$ is determined by the
periodicity condition $c[s,s+n-2]=c[t,t+n-2]$. This fact implies that
$U_2(G_n) \geq n-1$, and, thus, we have that $U(G_n) \geq n-1$, which,
in particular, gives condition (C4).

Let us check that $R(G_n) \leq n + R(G_1)$, which will imply that
$(G_n)$ satisfies condition (C5). The statement that $R(G_n) \leq n +
R(G_1)$ is equivalent to the statement that for any two paths $b,c \in
B_{n-1}(G_1)$, there exists a path $d$ in $G_1$ of length less than or
equal to $R(G_1)$ such that $bdc$ is a path in $G_1$. In this
formulation, the statement is clearly true, since, by the definition
of~$R(G_1)$, there is a path $d$ from $t(b)$ to $i(c)$ of length less than
or equal to~$R(G_1)$, and then the concatenation $bdc$ gives a path in $G_1$.

Let $w^1$ be a positive left (row) eigenvector for $A_1$ (corresponding
to the eigenvalue $\lambda$), and let $v^1$ be a positive right
(column) eigenvector for~$A_1$ (corresponding to the eigenvalue $\lambda
$). Let $b \in B_{n-1}(G_1) = V_n$. Then let $w^n_b = w^1_{i(b)}$ and
$v^n_b = v^1_{t(b)} \lambda^{-(n-1)}$. Then $w^n$ is a positive left
eigenvector for $A_n$ and $v^n$ is a positive right eigenvector for
$A_n$. It follows that $(G_n)$ \mbox{satisfies}~condi\-tions~(C6) and (C7). In
fact, to satisfy (C7), we may choose $K = \max(K_1, K_2)$, where $K_1 =
(\max_i w^1_i) (\min_i w^1_i)^{-1}$ and $K_2 = (\max_i v^1_i) (\min_i
v^1_i)^{-1}$.

Condition (C8) follows from the fact that $(G_n)$ satisfies condition
(C7) (by applying Lemmas \ref{C7impliesC1C6} and \ref{UnifExpLemma} in
succession).
\end{pf}

\subsection{Probabilistic framework} \label{ProbFramework}

Let $\Omega$ be the probability space consisting of the set $\{0,1\}^n$
and the probability measure $\PPP$, where $\PPP$ is the product of the
Bernoulli measures on each coordinate with parameter $\al\in[0,1]$.
There is a natural partial order on $\Omega$, given by the relation
$\omega\leq\tau$ if and only if $\omega_i \leq\tau_i$ for $i = 1,\ldots,
n$. We say that a random variable $\chi$ on $\Omega$ is \textit
{monotone increasing} if $\chi(\omega) \leq\chi(\tau)$ whenever
$\omega\leq\tau$. An event $A$ is monotone increasing if its
characteristic function is monotone increasing. Monotone decreasing is
defined analogously. Monotone random variables and events have been
studied extensively \cite{Gr2}; however, we require only a small
portion of that theory. In particular, we will make use of the
following proposition, a proof of which may be found in \cite{Gr2}.
\begin{prop}[(FKG Inequality)] \label{FKG}
If $X$ and $Y$ are monotone increasing random variables on $\{0,1\}^n$,
then $\E(XY) \geq\E(X)\E(Y)$.
\end{prop}

It follows easily from the FKG Inequality that if $\bigcap F_j$ is a
finite intersection of monotone decreasing events, then $\PPP(\bigcap F_j)
\geq\prod\PPP(F_j)$ (use induction and note that if $\chi_F$ is the
characteristic function of the monotone decreasing event~$F$, then
$-\chi_F$ is monotone increasing). In fact, we only use this corollary,
but we nonetheless refer to it as the FKG Inequality.

For a finite, directed graph $G$, we consider the discrete probability
space on the set $\Omega_G = \{0,1\}^E$, where $\PPP$ is the product of
the Bernoulli($\al$) measures on each coordinate. The set $\Omega_G$
corresponds to the power set of $E$ in the usual way: $\omega$ in
$\Omega_G$ corresponds to the set $F$ in $2^E$ such that $e$ is in $F$
if and only if $\omega(e) = 1$. Furthermore, $\Omega_G$ corresponds to
the space of subgraphs of $G$: for $\omega$ in $\Omega_G$, define the
subgraph $G(\omega)$ to have vertex set~$V$ and edge set $F_{\omega}$,
where an edge $e$ in $E$ is included in $F_{\omega} \subset E$ if and
only if $\omega(e)=1$. In the percolation literature, the edges $e$
such that $\omega(e)=1$ are often called ``open,'' and the remaining
edges are called ``closed.'' Since we are interested in studying edge
shifts defined by graphs, we will refer to an edge~$e$ as ``allowed''
when $\omega(e)=1$ and ``forbidden'' when $\omega(e)=0$. Finally,
each~$\omega$ in $\Omega_G$ can be associated to the SFT $X_{\omega}$
defined as the set of all bi-infinite, directed walks on $G$ that
traverse only \textit{allowed} edges (with respect to $\omega$). The
probability measure $\PPP$ corresponds to allowing each edge of $G$
with probability $\al$, independently of all other edges. For the sake
of notation, we suppress the dependence of $\PPP$ on the graph $G$.
%
\begin{defn}
In this work we consider the following conjugacy invariants of SFTs.
Let~$\cE$ be the property containing only the empty shift. Let~$\cZ$ be
the property containing all SFTs with zero entropy. By convention, we
let $\cE\subset\cZ$. For any SFT $X$, let $\h(X)$ be the topological
entropy, and let $I(X)$ be the number of irreducible components of $X$.
If $X$ is nonempty, let $\beta(X)$ be defined by the equation $\h(X) =
\log(\beta(X))$. If $X$ is empty, let $\beta(X) = 0$. If $\cS$ is a
property of SFTs and $G$ is a finite directed graph, then let $\cS_G
\subset\Omega_G$ be the set of $\omega$ in $\Omega_G$ such that
$X_{\omega}$ has property $\cS$. If $f$ is a~function from SFTs to the
real numbers and $G$ is a finite directed graph, then let $f_G\dvtx \Omega
_G \to\R$ be the function $f_G(\omega) = f(X_{\omega})$.
\end{defn}

\section{Emptiness} \label{Emptiness}

Recall that $\Sp_{\times}(G)$, $\zeta_G$ and $z(G)$ were defined in
Definitions \ref{nonzeroSpec}, \ref{zetaDefn} and \ref{Defnofz}, respectively.
\begin{theorem} \label{EmptyThm}
Let $(G_n)$ be a sequence of graphs such that $\Sp_{\times}(G_n) = \Sp
_{\times}(G_1)$ for all $n$ and either \textup{(i)} $\lambda= \lambda_{G_1} =
1$ or \textup{(ii)} $\lambda= \lambda_{G_1} > 1$ and $z(G_n)$ tends to infinity
as $n$ tends to infinity. Let $\zeta= \zeta_{G_1}$. Then
\[
\lim_{n \rightarrow\infty} \PPP(\mathcal{E}_{G_n}) = \cases{(\zeta(\al
))^{-1}, &\quad if $\al\in[0, 1/\lambda)$,\cr
0, &\quad if $\al\in[1/\lambda,1]$.}
\]
\end{theorem}
\begin{rmk} \label{EmptyRmk} Theorem~\ref{EmptyThmIntro} can be
obtained as a corollary of Theorem~\ref{EmptyThm} by taking $(G_n)$ to
be the sequence of $n$-block graphs of $X$. Indeed, if the SFT $X$ in
Theorem \ref{EmptyThmIntro} has zero entropy, then $\lambda= 1$, and
the conclusion of Theorem \ref{EmptyThmIntro} follows from case (i) in
Theorem~\ref{EmptyThm}. If the SFT $X$ in Theorem~\ref{EmptyThmIntro}
has positive entropy, then $\lambda>1$ and $z(G_n)$ tends to infinity
by the exact same argument in the proof of Proposition \ref
{higherGraphSeqLemma}(iii), and, therefore, the conclusion of Theorem
\ref{EmptyThmIntro} follows from case (ii) in Theorem \ref{EmptyThm}.
\end{rmk}

In this section we provide a proof of Theorem \ref{EmptyThm}. Before
proceeding with the proof, we state a fact that will be useful in the
investigations that follow. Recall that for a path $b$, we denote by
$V(b)$ the set of vertices traversed by~$b$.
\begin{lemma} \label{perLemma}
Suppose $G$ is a directed graph. Suppose $b$ is in $\Per(G)$ such that
$|V(b)| < \per(b)$. Then there exists a path $c$ in $\Per(G)$ such that
$\per(c) < \per(b)$ and $V(c) \subset V(b)$.
\end{lemma}
\begin{pf}
Let $v$ be in $V(b)$. Then there exists a return path to $v$ following~%
$b$, and we may choose a shortest return path $c$ to $v$ using only
vertices in~$V(b)$. Then $c$ is in $\Per(G)$ and $\per(c) < \per(b)$,
as desired.
\end{pf}
\begin{pf*}{Proof of Theorem \ref{EmptyThm}}
Recall that an SFT is nonempty if and only if it contains a periodic
point (see \cite{LM}).

First, assume that case (i) holds, which means that $\lambda= 1$. In
this case, each $X_{G_n}$ contains finitely many orbits. Further, the
number of periodic orbits of each period in $X_{G_n}$ is constant, and
the probability of each periodic orbit being allowed in $X_{\omega}$ is
constant. Therefore, the conclusion follows immediately, since the
sequence $\PPP(\mathcal{E}_{G_n})$ is constant.

Now assume that case (ii) holds. For the moment, consider a fixed
natural number $n$. Let $\{\gamma_j \}_{j\in\N}$ be an enumeration of
the periodic orbits of $X_{G_n}$ such that if $i \leq j$, then $\per
(\gamma_i) \leq\per(\gamma_j)$. Let $p_i = \per(\gamma_i) = |\gamma
_i|$. Let $V_n(\gamma_j)$ be the vertices in $G_n$ traversed in the
orbit $\gamma_j$ and let $E_n(\gamma_j)$ be the edges in~$G_n$
traversed in the orbit $\gamma_j$.

Now for each $j$, let $A_j$ be the event that $\gamma_j$ is allowed,
which is the event that all of the edges in $E_n(\gamma_j)$ are
allowed. Let $F_j$ be the event that $\gamma_j$ is forbidden, which is
$A_j^c$, the complement of $A_j$. Notice that $A_j$ is a monotone
increasing event (if $\omega$ is in $A_j$ and $\omega\leq\omega'$,
then $\omega'$ is in $A_j$), and $F_j$ is a~monotone decreasing event.
The fact that an SFT is nonempty if and only if it contains a periodic
point implies that $\cE_{G_n} = \bigcap F_j$.

Combining the definition of $z(G_n)$ and Lemma \ref{perLemma}, we
obtain that if $\per(\gamma_i) \leq z(G_n)$, then $|E_n(\gamma_i)| =
p_i$. It follows that $\PPP(F_i) = 1-\al^{p_i}$ for each~$i$ such that
$p_i \leq z(G_n)$. Furthermore, the definition of $z(G_n)$ implies that
the events $F_i$ such that $p_i \leq z(G_n)$ are all jointly
independent. These observations give that
%
\begin{eqnarray}
\label{EmptyUpperBd}
\PPP(\cE_{G_n}) &=& \PPP\biggl(\bigcap_{j \in\N} F_j\biggr) \leq\PPP\biggl( \bigcap
_{p_i \leq z(G_n)} F_i \biggr) \\
\label{EmptyUpperBd2}
&=& \prod_{p_i \leq z(G_n)} \PPP(F_i) =
\prod_{p_i \leq z(G_n)} (1-\al^{p_i}).
\end{eqnarray}

Using Lemma \ref{perLemma}, we see that there is great redundancy in
the intersection~$\bigcap F_j$. Eliminating some of this redundancy, we
obtain the following:
%
\begin{equation} \label{infiniteCapToFinite}
\bigcap_{j \in\N} F_j = \bigcap_{j\dvtx |E_n(\gamma_j)|=p_j} F_j.
\end{equation}
Then using Lemma \ref{perLemma} again and the fact that $|E_n(\gamma
_j)| \leq|E_n|$, we see that the intersection\vadjust{\goodbreak} on the right in
(\ref{infiniteCapToFinite}) is actually a finite intersection. Applying
the FKG Inequality, we obtain that
%
\begin{eqnarray}
\label{EmptyLowerBd1}\qquad
\PPP(\cE_{G_n}) & = & \PPP\biggl(\bigcap_{j\in\N} F_j\biggr)
= \PPP\biggl( \bigcap_{j\dvtx |E_n(\gamma_j)|=p_j} F_i \biggr) \geq\prod_{j\dvtx |E_n(\gamma_j)|=p_j} \PPP
(F_i) \\
\label{EmptyLowerBd2}
& = & \prod_{j\dvtx |E_n(\gamma_j)|=p_j}
(1-\al^{p_j}) \geq\prod_{j\dvtx p_j \leq|E_n|} (1-\al^{p_j}) .
\end{eqnarray}

Combining the inequalities in (\ref{EmptyUpperBd}), (\ref
{EmptyUpperBd2}), (\ref{EmptyLowerBd1}) and (\ref{EmptyLowerBd2}) gives
that for each $n$,
%
\begin{equation} \label{EmptyCombinedBds}
\prod_{ p_j \leq|E_n|} (1-\al^{p_j}) \leq\PPP(\cE_{G_n}) \leq\prod
_{p_i \leq z(G_n)} (1-\al^{p_i}).
\end{equation}

By the standing assumptions that $\Sp_{\times}(G_n) = \Sp_{\times
}(G_1)$, we have that $|{\Per_p}(G_n)\hspace*{-0.4pt}|$ is independent of $n$. Since
$z(G_n)$ and $|E_n|$ tend to infinity as~$n$ tends to infinity,
equation (\ref{EmptyCombinedBds}) gives that
\[
\lim_{n \to\infty} \PPP(\cE_{G_n}) = \prod_{j=1}^{\infty} (1-\al^{p_j}).
\]
Then Theorem \ref{EmptyThm} follows from the well-known product formula
for $\zeta$ (see~\cite{LM}), which may be stated as
\[
(\zeta(t))^{-1} = \prod_{j=1}^{\infty} (1-t^{p_j}),
\]
along with the fact that $\zeta(t)$ converges for $t < 1/\lambda$ and
diverges to $+\infty$ for $t \geq1/\lambda$.
\end{pf*}

\section{Subcritical phase} \label{subcriticalPhase}

In this section we study random SFTs in the subcri\-tical phase: $0 \leq
\al< 1/\lambda$. The main result of this section is Theorem~\ref
{SubCritThm}. Let us fix some notation for this section. We consider a
sequence of graphs~$(G_n)$ such that $\Sp_{\times}(G_n) = \Sp_{\times
}(G_1)$ and $z(G_n)$ tends to infinity as $n$ tends to infinity, with
$\lambda= \lambda_{G_1} \geq1$ and $\zeta= \zeta_{G_1}$. Since $\Sp
_{\times}(G_n) = \Sp_{\times}(G_1)$, there exist shift-commuting
bijections $\phi_n\dvtx \Per(X_{G_1}) \to\Per(X_{G_n})$. In other words,
there exist bijections $\phi_n$ from the set of cyclic paths in $G_1$
to the set of cyclic paths in $G_n$ such that if $b$ is in $\Per
_p(G_1)$, then $\phi_n(b)$ is $\Per_p(G_n)$. If $b$ is in $\Per(G)$,
then we refer to $\theta(b)$ (recall Definition \ref{PerDefn}) as a
cycle. Using the fixed bijections $\phi_n$, we may refer to a cycle
$\gamma$ as being in $G_n$ for any $n$. We fix an enumeration of the
cycles in $G_1$, $\{\gamma_i\}_{i \in\N}$, and then since the
bijections~$\phi_n$ are fixed, this choice simultaneously gives
enumerations of all the cycles in each $G_n$. For any $s$ in $\N$, let
$p_s = \per(\gamma_s)$. Let us begin with a lemma.
\begin{lemma} \label{finiteIrrCompLemma}
Let $(G_n)$ be a sequence of graphs such that $\Sp_{\times}(G_n) = \Sp
_{\times}(G_1)$ and $z(G_n)$ tends to infinity as $n$ tends to
infinity, with $\lambda= \lambda_{G_1} \geq1$ and\vadjust{\goodbreak} $\zeta= \zeta
_{G_1}$. Given a nonempty, finite set $S$ in $\N$, let $D_{G_n}(S)$ be
the event that the set of allowed cycles is $\{ \gamma_s\dvtx s \in S\}$. Then
\[
\lim_{n \rightarrow\infty} \PPP(D_{G_n}(S)) = \cases{\displaystyle (\zeta(\al))^{-1}
\prod_{j \in S} \frac{\al^{p_j}}{1-\al^{p_j}}, &\quad
if $\al\in[0, 1/\lambda)$, \vspace*{3pt}\cr
0, &\quad if $\al\in[1/\lambda,1]$.}
\]
\end{lemma}

The proof of Lemma \ref{finiteIrrCompLemma} is an easy adaptation of
the proof of Theorem~\ref{EmptyThm}, and we omit it for the sake of brevity.

Recall that $I(X)$ denotes the number of irreducible components in the
SFT $X$, and for any graph $G$, the random variable $I_G\dvtx \Omega_G \to
\Z_{\geq0}$ is defined by the equation $I_{G}(\omega) = I(X_{\omega})$.
\begin{theorem} \label{SubCritThm}
Let $(G_n)$ be a sequence of graphs such that $\Sp_{\times}(G_n) = \Sp
_{\times}(G_1)$ and either \textup{(i)} $\lambda= \lambda_{G_1} = 1$ or
\textup{(ii)}
$\lambda= \lambda_{G_1} > 1$ and $z(G_n)$ tends to infinity as $n$
tends to infinity. Let $\zeta= \zeta_{G_1}$. Then for $0\leq\al<
1/\lambda$,
\begin{longlist}[(2)]
\item[(1)] $\lim_{n \rightarrow\infty} \PPP(\mathcal{Z}_{G_n}) = 1$;
\item[(2)] the sequence $(I_{G_n})$ converges in distribution to the random
variable~$I_{\infty}$ such that $\mathbb{P}(I_{\infty} = 0) = (\zeta(\al
))^{-1}$ and for $k \geq1$,
\[
\mathbb{P}(I_{\infty} = k) = (\zeta(\al))^{-1} \mathop{\sum_{S \subset
\N}}_{|S|=k} \prod_{s \in S} \frac{ \al^{p_s}}{1-\al^{p_j}},
\]
where $\{\gamma_i \}_{i=1}^{\infty}$ is an enumeration of the cycles in $G_1$;
\item[(3)] the random variable $I_{\infty}$ has exponentially decreasing
tail and therefore finite moments of all orders.
\end{longlist}
\end{theorem}
\begin{rmk} One obtains Theorem~\ref{SubCritThmIntro} as a consequence
of Theorem~\ref{SubCritThm} by taking $(G_n)$ to be the sequence of
$n$-block graphs of a nonempty SFT $X$. Indeed, if the SFT $X$ in
Theorem \ref{SubCritThmIntro} has zero entropy, then $\lambda= 1$, and
the conclusions of Theorem \ref{SubCritThmIntro} follow from the case
(i) in Theorem \ref{SubCritThm}. If the SFT $X$ in
Theorem~\ref{SubCritThmIntro} has positive entropy, then $\lambda>1$ and $z(G_n)$
tends to infinity by the exact same argument in the proof of
Proposition~\ref{higherGraphSeqLemma}(iii), and, therefore, the
conclusions of Theorem \ref{SubCritThmIntro} follow from case (ii) in
Theorem~\ref{SubCritThm}.
\end{rmk}

\begin{pf*}{Proof of Theorem \ref{SubCritThm}}
Let $(G_n)$ be as above. Let $0 \leq\al< 1/\lambda$.

First, assume that case (i) holds, which means that $\lambda= 1$.
Conclusion (1) follows immediately, since for each $n$, we have that
$\PPP(\mathcal{Z}_{G_n}) = 1$ [the random SFT $X_{\omega}$ satisfies $0
= \h(X_{\omega}) \leq\h(X_{G_n}) = \log\lambda= 0$]. Also, the fact
that $\lambda= 1$ is equivalent to the fact that $G_1$ (and therefore
$G_n$) contains only finitely many cycles. Then conclusions (2) and (3)
also follow immediately, since the sequence $I_{G_n}$ is constant.

Now assume that case (ii) holds. Recall that we have an enumeration $\{
\gamma_i\}_{i \in\N}$ of the cycles in $G_1$,\vadjust{\goodbreak} which we refer to as an
enumeration of the cycles in $G_n$, for any~$n$, using the bijections
$\phi_n$. Also recall that for any nonempty, finite set $S \subset\N$,
we denote by $D_{G_n}(S)$ the event in $\Omega_{G_n}$ consisting of all
$\omega$ such that the set of cycles in $G_n(\omega)$ is exactly $\{
\gamma_s\dvtx s\in S\}$.

\textit{Proof of} (1).
Recall that an SFT has zero entropy if and only if it has at most
finitely many periodic points \cite{LM}. Then we have that
%
\begin{equation} \label{ZeroEntUnion}
\cZ_{G_n} = \cE_{G_n} \cup\biggl( \mathop{\bigcup_{S \subset\N}}_{0 < |S| <
\infty} D_{G_n}(S) \biggr).
\end{equation}
Also note that by the definition of $D_{G_n}(S)$, the union in (\ref
{ZeroEntUnion}) is a disjoint union. Thus, we have that
\[
\PPP(\cZ_{G_n}) = \PPP(\cE_{G_n}) + \mathop{\sum_{S \subset\N}}_{0 <
|S| < \infty} \PPP(D_{G_n}(S)).
\]
Now let $S_1,\ldots, S_J$ be distinct, nonempty, finite subsets of $\N$.
Then by Theorem~\ref{EmptyThm} and Lemma \ref{finiteIrrCompLemma}, we
have that
\begin{eqnarray*}
\liminf_{n \to\infty} \PPP(\cZ_{G_n}) & \geq &\lim_{n \to\infty} \PPP
(\cE_{G_n}) + \sum_{j=1}^J \lim_{n \to\infty} \PPP(D_{G_n}(S_j)) \\
& = & (\zeta(\al))^{-1}\Biggl(1 + \sum_{j=1}^J \prod_{s \in S_j} \frac{\al
^{p_s}}{1-\al^{p_s}}\Biggr).
\end{eqnarray*}
Since $J$ and $S_1,\ldots, S_J$ were arbitrary, we conclude that
\[
\liminf_{n \to\infty} \PPP(\cZ_{G_n}) \geq(\zeta(\al))^{-1}\biggl( 1+
\mathop{\sum_{S \subset\N}}_{0 < |S| < \infty} \prod_{s \in S} \frac
{\al^{p_s}}{1-\al^{p_s}} \biggr).
\]
Using the facts that $\al^{p_s}/(1-\al^{p_s}) = \sum_{k=1}^\infty(\al
^{p_s})^k$ and $\al< 1/\lambda$ (which implies that the relevant
infinite products and series converge uniformly), one may easily check that
\[
\biggl( 1+ \mathop{\sum_{S \subset\N}}_{0 < |S| < \infty} \prod_{s \in S}
\frac{\al^{p_s}}{1-\al^{p_s}} \biggr) = \zeta(\al).
\]
Thus, we have shown that $\liminf_{n} \PPP(\cZ_{G_n}) \geq1$. Since
$\limsup_{n} \PPP(\cZ_{G_n}) \leq1$, we conclude that $\lim_{n} \PPP
(\cZ_{G_n}) = 1$.

\textit{Proof of} (2).
Since $I_{G_n}$ takes values in $\Z_{\geq0}$, the sequence $(I_{G_n})$
converges in distribution to $I_{\infty}$ if and only if $\PPP
(I_{G_n}=k)$ converges to $\PPP(I_{\infty}=k)$ for each $k$ in $\Z
_{\geq0}$.

Note that $I_{G_n}(\omega) = 0$ if and only if $\omega$ is in $\cE
_{G_n}$, which implies that $\PPP(I_{G_n}=0) = \PPP(\cE_{G_n})$. Thus,
for $\al< 1/\lambda$, Theorem \ref{EmptyThm} implies that\break $\PPP
(I_{G_n}=0)$ converges to $(\zeta(\al))^{-1}$ as $n$ tends to
infinity.\vadjust{\goodbreak}

Now let $k$ be in $\N$. Recall that $\{\gamma_i\}_{i=1}^{\infty}$ is an
enumeration of the cycles in $G_1$, and we have fixed bijections
between these cycles and the cycles in each $G_n$. By Theorem \ref
{SubCritThm}(1), we have that $\lim_n \PPP(\cZ_{G_n}) =1$, and,
therefore, $ \PPP(I_{G_n} = k) = \PPP(\{I_{G_n} = k\} \cap\cZ_{G_n}) +
\varepsilon_n$, where $\varepsilon_n$ tends to $0$ as $n$ tends to
infinity. Thus, we need only focus on events of the form $\{I_{G_n} =
k\} \cap\cZ_{G_n}$ for some $k$. Now if $\omega$ is in~$\cZ_{G_n}$,
then $I_{G_n}(\omega)$ is the number of periodic orbits in~$X_{\omega
}$. Thus,
\[
\PPP(\{I_{G_n} = k\} \cap\cZ_{G_n}) = \mathop{\sum_{S \subset\N
}}_{|S|=k} \PPP(D_{G_n}(S)).
\]

For any $n$ in $\N$, let $T^0_n = \PPP(\cE_{G_n})$. For $k$ in $\N$ and
$n$ in $\N$, let
\[
T^k_n = \mathop{\sum_{S \subset\N}}_{|S|=k} \PPP(D_{G_n}(S)).
\]
We have\vspace*{2pt} that $\sum_{k=0}^{\infty} T_n^k = \PPP(\cZ_{G_n})$, and,
therefore, $\lim_n \sum_{k=0}^{\infty} T_n^k = 1$ by Theorem \ref
{SubCritThm}(1). Also, using Lemma \ref{finiteIrrCompLemma}, we have
that $\liminf_n T^k_n \geq T^k$, where $T^0 = (\zeta(\al))^{-1}$ and
for $k$ in $\N$,
\[
T^k = (\zeta(\al))^{-1} \mathop{\sum_{S \subset\N}}_{|S| = k} \prod_{s
\in S} \frac{\al^{p_s}}{1-\al^{p_s}}.
\]
Further, we have that $\sum_{k=0}^{\infty} T^k = 1$. It follows from
these facts that $\lim_n T^k_n = T^k$. Thus, we have shown that for $k$
in $\N$,
\[
\lim_n \PPP( I_{G_n} = k) = \lim_n \PPP(\{I_{G_n} = k\} \cap\cZ_{G_n})
= (\zeta(\al))^{-1}
\mathop{\sum_{S \subset\N}}_{|S| = k} \prod_{s \in S} \frac{\al
^{p_s}}{1-\al^{p_s}}
\]
as desired.

\textit{Proof of} (3).
For $k$ in $\N$, let
\[
T^k = \PPP(I_{\infty} = k) = (\zeta(\al))^{-1} \mathop{\sum_{S \subset
\N}}_{|S| = k} \prod_{s \in S} \frac{\al^{p_s}}{1-\al^{p_s}}.
\]
We show that there for any real number $\delta>0$, there exists $k_0$
such that $T^{k+1} \leq\delta T^k$ for all $k \geq k_0$. Let $\delta
>0$. Since $\al< 1/\lambda$, we have that
\[
\sum_{i \in\N} \frac{\al^{p_i}}{1-\al^{p_i}} < \infty.
\]
Now choose $k_0$ such that
\[
\sum_{i \geq k_0} \frac{\al^{p_i}}{1-\al^{p_i}} < \delta.
\]
In the following sums, we will use that any set $S \subset\N$ with
$|S|=j$ can be written as $S = \{s_1,\ldots, s_j\}$, where $s_1 < \cdots
< s_j$. Note that in this case $s_j \geq j$. Then for $k \geq k_0$ we have
\begin{eqnarray*}
(\zeta(\al)) T^{k+1} & = & \mathop{\sum_{S \subset\N}}_{|S| = k+1} \prod
_{i=1}^{k+1} \frac{\al^{p_{s_i}}}{1-\al^{p_{s_i}}}
= \mathop{\sum_{S \subset\N}}_{|S| = k} \prod_{i=1}^{k} \frac{\al
^{p_{s_i}}}{1-\al^{p_{s_i}}} \sum_{j > s_k} \frac{\al^{p_{j}}}{1-\al
^{p_{j}}} \\
& \leq & \mathop{\sum_{S \subset\N}}_{|S| = k} \prod_{i=1}^{k} \frac{\al
^{p_{s_i}}}{1-\al^{p_{s_i}}} \sum_{j > k_0} \frac{\al^{p_{j}}}{1-\al^{p_{j}}}
\leq \Biggl(\mathop{\sum_{S \subset\N}}_{|S| = k} \prod_{i=1}^{k} \frac{\al
^{p_{s_i}}}{1-\al^{p_{s_i}}} \Biggr) \delta \\
& = & (\zeta(\al)) T^k \delta.
\end{eqnarray*}
Since $\al< 1/\lambda$, we have that $0 < \zeta(\al) < \infty$, and we
conclude that $T^{k+1} \leq\delta T^k$ for all $k \geq k_0$.
\end{pf*}

We recognize the distribution of $I_{\infty}$ as the sum of countably
many independent Bernoulli trials, where the probability of success of
trial $i \in\N$ is given by $\al^{p_i}$ for some enumeration $\{\gamma
_i\}_{i\in\N}$ of the cycles in $G_1$ (or any $G_n$). We record some
facts about this distribution in the following corollary.
\begin{cor}
With the same hypotheses as in Theorem \ref{SubCritThm}, the
characteristic function of $I_{\infty}$ is given by
\[
\varphi_{I_{\infty}}(t) = (\zeta(\al))^{-1} \prod_s \biggl( 1 + e^{it} \frac{
\al^{p_s}}{1-\al^{p_s}} \biggr),
\]
where the product is over all periodic orbits in $X$. It follows that
the moment generating function of $I_{\infty}$ is given by
\[
M_{I_{\infty}}(t) = (\zeta(\al))^{-1} \prod_s \biggl( 1 + e^{t} \frac{ \al
^{p_s}}{1-\al^{p_s}} \biggr).
\]
\end{cor}
\begin{rmk}
In Theorems \ref{EmptyThm} and \ref{SubCritThm}, we assert the
existence of various limits to certain values. Beyond the bounds given
in our proofs, we do not know at which rates these sequences converge
to their limits.
\end{rmk}

\section{Supercritical phase} \label{supercriticalPhase}

In this section we study random SFTs in the supercritical phase. The
main results are Theorems \ref{EntropyThm} and \ref{MSFTthm}. On
a~first reading, the reader may prefer to skip Section \ref{ITSection}
and refer back to it as necessary. Our proof of Theorem \ref
{EntropyThm} relies, in part, on showing that with large probability
the number of allowed words of length $k$ in a random SFT is close to
$(\al\lambda)^k$, for a particular choice of $k$. In our proof, we
choose~$k$ to be polynomial in $m = m(G_n)$\vadjust{\goodbreak} for two reasons. First, we
need $k$ to dominate~$m$, so that the $k$th root of the number of words
of length $k$ gives a~good upper bound on the Perron eigenvalue of the
random SFT. Second,~$k$ should be subexponential in $m$, essentially
because most paths in $G_n$ with length subexponential in $m$ are
self-avoiding, and we need good bounds on the probability of paths of
length $k$ that exhibit ``too-soon-recurrence.'' For context, we recall
a result of Ornstein and Weiss \cite{OW}. In fact, their result is
quite general, but we only recall it in a very specific case. Let $X$
be an irreducible SFT with measure of maximal entropy $\mu$. For $x$ in
$X$, let $R_n(x)$ be the first return time (greater than 0) of $x$ to
the cylinder set $x[1,n]$ under~$\sigma$. Then the result of Ornstein
and Weiss implies that for $\mu$-a.e. $x$ in $X$, $\lim_n n^{-1} \log
R_n(x) = h(X)$. It follows from this result that for $k$ polynomial in
$n$, the $\mu$-measure of the set of words of length $k$ with a
repeated $n$-word tends to $0$. In the following lemmas, we give some
quantitative bounds on the $\mu$-measure of the set of paths of length
$k$ in $G_n$ with $k-j$ repeated edges, where the important point for
our purposes is that the bounds improve exponentially as $j$ decreases.
To get these bounds, we employ some of the language and tools of
information theory. After getting a handle on the $\mu$-measure of
paths in $G_n$ with certain self-intersection properties, our
assumption that $(G_n)$ satisfies condition (C7) in Definition \ref{conditions}
implies that the $\mu$-measure on paths is the same as the counting
measure up to uniform constants.

\subsection{Information theory and lemmas} \label{ITSection} In keeping
with the convention of information theory, $\log(x)$ denotes the base
$2$ logarithm of $x$.
\begin{defn}
A binary $n$-\textit{code} on an alphabet $\cA$ is a mapping $C\dvtx\allowbreak \cA
^{n} \to\{0,1\}^*$, where $\{0,1\}^*$ is the set of all finite words
on the alphabet $\{0,1\}$. We may refer to such mappings simply as
codes. A code is \textit{faithful} if it is injective. The function
that assigns to each $w$ in $\cA^{n}$ the length of the word~$C(w)$ is
called the \textit{length function} of the code, and it will be denoted
by $\mathcal{L}$ when the code is understood. A~code is a \textit
{prefix code} if \mbox{$w = w'$} whenever~$C(w)$ is a prefix of $C(w')$. A
\textit{Shannon code} with respect to a~measure~$\nu$ on~$\cA^n$ is
a~code such that $\mathcal{L}(w) = \lceil- \log\nu(w) \rceil$.
\end{defn}

We note that for a measure $\nu$ on $\cA^n$, there is a prefix Shannon
code on~$\cA^n$ with respect to $\nu$ \cite{Shields}. We will also
require the following two lemmas from information theory.
\begin{lemma}[\cite{Shields}] \label{mainInfoTheoryLemma}
Let $\cA$ be an alphabet. Let $C_n$ be a prefix-code on $\cA^n$, and
let $\mu$ be a shift-invariant Borel probability measure on $\cA^{\Z}$. Then
\[
\mu\bigl( \{ w \in\cA^n\dvtx \mathcal{L}(w) + \log\mu( w ) \leq-a\} \bigr) \leq2^{-a}.
\]
\end{lemma}
\begin{pf}
Let $B = \{ w \in\cA^n\dvtx \mathcal{L}(w) + \log\mu( w ) \leq-a\}$.
Then for any $w$ in $B$, we have that $\mu(w) \leq2^{-\mathcal{L}(w)}
2^{-a}$. The Kraft inequality\vadjust{\goodbreak} for prefix codes (\cite{Shields}, page 73)
states that since $\mathcal{L}$ is a prefix code, $\sum_{w \in\cA^n}
2^{-\mathcal{L}(w)} \leq1$. Hence,
\[
\mu(B) = \sum_{w \in B} \mu(w) \leq2^{-a} \sum_{w \in B} 2^{-\mathcal
{L}(w)} \leq2^{-a}.
\]
\upqed\end{pf}
\begin{lemma}[\cite{Shields}] \label{EliasLemma}
There is a prefix code $C\dvtx \N\to\{0,1\}^*$ such that $\ell(C(n)) =
\log(n) + o(\log(n))$, where $\ell(C(n))$ is the length of $C(n)$.
\end{lemma}
\begin{defn}
A prefix code satisfying the conclusion of Lemma \ref{EliasLemma} is
called an \textit{Elias code}.
\end{defn}

Recall that if $b$ is a path in the graph $G = (V,E)$, then we denote
by $E(b)$ the set of edges traversed by $b$. Let $(G_n)$ be a sequence
of graphs satisfying our Standing Assumptions \ref{StAssump}.
\begin{defn} \label{NDef}
For each $n$, $k$, and $1 \leq j \leq k-1$, let
\[
N_{n,k}^j = \{ b \in B_k(G_n)\dvtx |E_n(b)|\leq j\}.
\]
\end{defn}
\begin{defn} \label{DDef}
For each $n$, $k$, and $1 \leq j \leq2k-1$, let
\[
D_{n,k}^j = \{ (b,c) \in B_k(G_n) \times B_k(G_n)\dvtx E_n(b) \cap E_n(c)
\neq\varnothing, |E_n(b) \cup E_n(c)| \leq j\}.
\]
\end{defn}
\begin{defn} \label{QDef}
For each $n$, $k$, and $1 \leq j \leq k-1$, let
\[
Q_{n,k}^j = \{ b \in\Per_k(G_n)\dvtx |E_n(b)|\leq j\}.
\]
\end{defn}
\begin{defn} \label{SDef}
For each $n$, $k$, and $1 \leq j \leq2k-1$, let
\[
S_{n,k}^j = \{ (b,c) \in\Per_k(G_n) \times\Per_k(G_n)\dvtx E_n(b) \cap
E_n(c) \neq\varnothing, |E_n(b) \cup E_n(c)| \leq j\}.
\]
\end{defn}

For any of the sets defined in Definitions \ref{NDef}--\ref{SDef}, we
use a ``hat'' to denote the set with ``$\leq$'' replaced by ``$=$'' in the
definition. For example,
\[
\hat{N}_{n,k}^j = \{ b \in B_k(G_n)\dvtx |E_n(b)| = j\}.
\]
The ``hat'' notation will only appear in the proof of Theorem \ref
{EntropyThm}. The following four lemmas find bounds on $|N_{n,k}^j|$,
$|D_{n,k}^j|$, $|S_{n,k}^{2k-1}|$ and $|S_{n,k}^j|$.\vspace*{1pt}

The following lemma bounds the $\mu$-measure (and therefore the
cardinality) of the set of paths of length $k$ in $G_n$ that traverse
at most $j<k$ edges. The proof relies on a general principle in
information theory (made precise by Lemma \ref{mainInfoTheoryLemma}): a
set of words that can be encoded ``too efficiently'' must have small
measure. In order to use this principle, we find an efficient encoding
of the paths of length $k$ in $G_n$ that traverse at most $j$ edges.
The basic observation behind the coding is trivial: a path of length
$k$ that only traverses $j < k$ edges must have $k-j$ repeated edges.
Therefore, instead of encoding each of the $k-j$ repeated edges
explicitly, we simply encode some combinatorial data that specifies
when ``repeats'' happen and when the corresponding edges are first traversed.
\begin{lemma} \label{manyRepeatsLemma}
Let $(G_n)$ be a sequence of graphs satisfying the Standing Assumptions
\ref{StAssump} and such that $(G_n)$ has local uniqueness of paths and
bounded distortion of weights [conditions \textup{(C4)} and \textup{(C7)} in Definition \ref
{conditions}]. Then there exists a polynomial $p_0(x)$ and $n_0$ such
that for each $n \geq n_0$, $k > U(G_n)$ and $1 \leq j \leq k-1$,
\[
\mu(N_{n,k}^j) \leq p_0(k)^{\min( k-j, k/U(G_n) )} \lambda^{-(m(G_n)+k-j)}
\]
and
\[
|N_{n,k}^j| \leq p_0(k)^{\min(k-j, k/U(G_n))} \lambda^j.
\]
\end{lemma}
\begin{pf}
Consider $(G_n)$, $n$, $k$ and $j$ as in the hypotheses. Let $m =
m(G_n)$ and $U = U(G_n)$. A path $b$ in $N_{n,k}^j$ from vertex $s$ to
vertex $t$ contributes $w_s^n v_t^n \lambda^{-k}$ to $\mu(N_{n,k}^j)$.
The condition (C7) gives a uniform constant $K$ such that $w_s^n v_t^n$
is bounded below by $(K^2 |V_n|)^{-1} = (K^2 \lambda^m)^{-1}$.
Therefore, the bound on $|N_{n,k}^j|$ follows from the bound on $\mu
(N_{n,k}^j)$, since $|N_{n,k}^j| \leq K^2 \lambda^{m+k} \mu(N_{n,k}^j)$
[as in Lem\-ma~\ref{implicationsLemma}(1)]. We now proceed to show the
bound on $\mu(N_{n,k}^j)$.

Let $r = k-j$. Consider $b$ in $N_{n,k}^j$. Then there exists $1<t_1 <
\cdots< t_r \leq k$ such that $b_{t_i} = b_{s_i}$ for some $1 \leq s_i
< t_i$, for each $i = 1,\ldots, r$, where $s_i = \min\{ s \geq1\dvtx b_s
= b_{t_i}\}$. Now we define a set $\I\subset\{1,\ldots, r\}$ by
induction. Let $i_1 = 1$ and $\I_1 = \{i_1\}$. Assuming by induction
that $i_j$ and $\I_j$ have been defined and that $i_j < r$, we define
$i_{j+1}$ and $\I_{j+1}$ as follows:
\begin{itemize}
\item if $t_{i_j+1} - t_{i_j} > U$, let $i_{j+1} = i_j +1$;
\item otherwise, if $t_{i_j+1} - t_{i_j} \leq U$, then let
\[
i_{j+1} = \max\{ i_j < i \leq r\dvtx t_i - t_{i_j} \leq U \}.
\]
\end{itemize}
Let $\I_{j+1} = \I_j \cup\{i_{j+1}\}$. This induction procedure
terminates when $i_j = r$ for some $j \leq r$, and we denote this
terminal $j$ by $j_*$. Let $\I= \I_{j_*}$. Note that for each $0 \leq
s \leq k-U$, we have that
\[
|\{ i \in\I\dvtx s+1 \leq t_i \leq s+U \} | \leq2.
\]
It follows that $|\I| \leq\min( r, 2 k/U + 2)$.

Having defined the set $\I$, we now decompose the integer interval $\{
1,\ldots, k\}$ into subintervals. First, let
\[
J = \bigcup_{j=1}^{j_*} \{t_{i_j}\} \cup\{ 1 \leq s \leq k\dvtx \exists i_j,
i_{j+1} \in\I, t_{i_{j+1}}-t_{i_j} \leq U \mbox{ and } t_{i_j} \leq s
\leq t_{i_{j+1}} \}.
\]
Let $J_1,\ldots, J_N$ be the maximal disjoint subintervals (with
singletons allowed) of $\{1,\ldots, k\}$ such that $J = J_1 \cup\cdots
\cup J_N$ and $J_{\ell} < J_{\ell+1}$. Note that $\sum_{\ell=1}^N
|J_{\ell}| = |J| \geq r$ and $N \leq|\I|$. Then let $I_1,\ldots,
I_{N+1}$ be the maximal disjoint subintervals of $\{1,\ldots, k\}$ such
that:
\begin{itemize}
\item$I_{\ell} \subset\{ 1,\ldots, k\} \setminus J$ for each $\ell=
1,\ldots, N+1$;\vspace*{2pt}
\item$\bigcup_{\ell=1}^{N+1} I_{\ell} = \{1,\ldots, k \} \setminus
J$;\vspace*{2pt}
\item and for each $\ell= 1,\ldots, N$, we have that $I_{\ell}$ is
nonempty and $I_{\ell} < I_{\ell+1}$.
\end{itemize}
In summary, we have that $\{1,\ldots, k\} = I_1 \cup J_1 \cup\cdots
\cup I_N \cup J_N \cup I_{N+1}$, and only $I_{N+1}$ may be empty.

For any $1 \leq s < t \leq k$, let $A_{s,t} = \{ b \in B_k(G_n)\dvtx b_s =
b_t\}$. By Lemma \ref{implicationsLemma}(3), there exists a uniform
constant $K_1$ such that
%
\begin{equation} \label{measOnPer}
\mu(A_{s,t}) \leq K_1 \lambda^{-m}.
\end{equation}

For notation, if $I$ is a subset of $\{1,\ldots,k\}$, then $b_{I}$ is
$b$ restricted to $I$. Since $\mu$ is a $1$-step Markov on $X_{G_n}$,
we have that
%
\begin{eqnarray}
\label{DecompOfMu1}
\mu(b |A_{s_1,t_1}) & = & \mu( b_{I_1} | A_{s_1,t_1}) \prod_{\ell=
1}^{N} \mu( b_{J_{\ell}} | A_{s_1,t_1} \cap b_{I_1 \ccdots I_{\ell}}
)\nonumber\\[-10pt]\\[-10pt]
&&{}\times\prod_{\ell= 2}^{N+1} \mu( b_{I_{\ell}} | A_{s_1,t_1} \cap b_{I_1
\ccdots J_{\ell-1}}) \nonumber\\[-2pt]
\label{DecompOfMu2}
& = & \mu( b_{I_1} | A_{s_1,t_1}) \prod_{\ell= 1}^{N} \mu( b_{J_{\ell}}
| b_{I_1 \ccdots I_{\ell}} ) \prod_{\ell= 2}^{N+1} \mu( b_{I_{\ell}} |
b_{I_1 \ccdots J_{\ell-1}}).
\end{eqnarray}

Given $b$, we may form $s_i$, $t_i$, $I_{\ell}$ and $J_{\ell}$ as
above, and then we encode $b$ as follows:
\begin{longlist}[(1)]
\item[(1)] encode $s_1$ and $t_1$ using an Elias code;
\item[(2)] encode $b_{I_1}$ using a prefix Shannon code with respect to $\mu
( \cdot| A_{s_1,t_1})$;
\item[(3)] assuming $b_{I_1 \ccdots I_{\ell}}$ has been encoded, we encode
$b_{J_{\ell}}$ by encoding $s_i$ and~$t_i$ for each $i$ in $\I$ such
that $t_{i} \in J_{\ell}$, using an Elias code (and note that this~%
infor\-mation completely determines $b_{J_{\ell}}$ by definition of $U$
and construction of~$J$);
\item[(4)] assuming $b_{I_1 \ccdots J_{\ell-1}}$ has been encoded, we encode
$b_{I_{\ell}}$ using a prefix Shannon code with respect to $\mu( \cdot
| b_{I_1 \ccdots J_{\ell-1}})$.
\end{longlist}

Now we analyze the performance of the code. Since the code is a
concatenation of prefix codes, it is a prefix code. Since $U$ tends to
infinity as $n$ tends to infinity [by (C4)] and $k > U$, there exists
$n_0$ such that for $n \geq n_0$ and $1 \leq s \leq k$, the length of
the codeword in the Elias encoding of $s$ is less than or equal to $2
\log k$. Then we have, neglecting bits needed to round up,
%
\begin{equation} \label{LengthOfCode}
\LL(b) \leq- \log\mu( b_{I_1} |A_{s_1,t_1}) + |\I|(4 \log k) + \sum
_{\ell=2}^{N+1} - \log\mu( b_{I_{\ell}} | b_{I_1 \ccdots
J_{\ell-1}}).\hspace*{-24pt}
\end{equation}
Combining (\ref{DecompOfMu1}), (\ref{DecompOfMu2}) and
(\ref{LengthOfCode}), we have that
%
\begin{equation}
\LL(b) + \log\mu(b) \leq|\I|(4 \log k) + \log\mu(A_{s_1,t_1}) + \sum
_{\ell=1}^N \log\mu( b_{J_{\ell}} | b_{I_1 \ccdots
I_{\ell}}).\hspace*{-25pt}\vadjust{\goodbreak}
\end{equation}
Now by Lemma \ref{implicationsLemma}(2) and (3), there exist uniform
constants $K_2$ and $K_3$ such that
%
\begin{eqnarray}
\LL(b) + \log\mu(b) & \leq & |\I|(4 \log k) + K_2 - m \log\lambda+ N
K_3 - |J| \log\lambda\\
& = & |\I|(4 \log k) + K_2 + N K_3 - (m + |J|) \log\lambda.
\end{eqnarray}
By construction, $|\I| \leq\min( k-j, 2k/U +2)$, $N \leq|\I|$ and
$|J| \geq r = k-j$. Then by Lemma \ref{mainInfoTheoryLemma}, there
exists a uniform constant $K_4>0$ such that
%
\begin{equation}
\mu(N_{n,k}^j) \leq(K_4 k^4)^{\min( k-j, 2k/U +2)} \lambda^{-(m+k-j)}.
\end{equation}
Letting $p_0(x) = K_5 x^{12}$, for some uniform constant $K_5>0$, we
obtain that
\[
\mu(N_{n,k}^j) \leq p_0(k)^{\min( k-j, k/U)} \lambda^{-(m+k-j)},
\]
which completes the proof.
\end{pf}

The following lemma bounds the $\mu\times\mu$-measure (and therefore
the cardinality) of the set of pairs paths of length $k$ in $G_n$ that
share at least one edge and together traverse at most $j<2k$ edges. The
general strategy of encoding pairs of paths using combinatorial data
and appealing to information theory is similar to that of Lemma \ref
{manyRepeatsLemma}. Lemma \ref{masterDoubleLemma} involves the
additional hypothesis that there exists a uniform bound $R$ such that
for any pair of paths $(u,w)$ in $G_n$, there exists a path $uvw$ in
$G_n$ with $|v|\leq R$. Using this hypothesis, one observes that pairs
of paths can essentially be concatenated in $G_n$ and then treated as
single paths as in Lem\-ma~\ref{manyRepeatsLemma}.
\begin{lemma} \label{masterDoubleLemma}
Let $(G_n)$ be a sequence of graphs satisfying the Standing Assumptions
\ref{StAssump} and such that $(G_n)$ has local uniqueness of paths,
small diameter and bounded distortion of weights [conditions \textup{(C4)},
\textup{(C5)} and \textup{(C7)} in Definition \ref{conditions}]. Then there exists a polynomial $p_1(x)$
and $n_1$ such that for $n \geq n_1$, $k > R(G_n)$ and $1 \leq j \leq2k-1$,
\[
\mu\times\mu(D^j_{n,k}) \leq p_1(k)^{\min(2k-j, k/U(G_n))} \lambda
^{-(m(G_n)+2k-j)}
\]
and
\[
|D^j_{n,k}| \leq p_1(k)^{\min(2k-j, k/U(G_n))} \lambda^{j+m(G_n)}.
\]
\end{lemma}
\begin{pf}
Consider $(G_n)$, $n$, $k$ and $j$ as in the hypotheses. Let $m =
m(G_n)$, $U=U(G_n)$ and $R = R(G_n)$. Note that the bound on
$|D_{n,k}^j|$ follows from the bound on $\mu\times\mu(D_{n,k}^j)$,
since condition (C7) implies that there exists a~uniform constant $K$
such that $|D_{n,k}^j| \leq K \lambda^{2m+2k} \mu\times\mu
(D_{n,k}^j)$ [as in Lem\-ma~\ref{implicationsLemma}(1)]. We now proceed
to show the bound on $\mu\times\mu(D_{n,k}^j)$.

By the definition of $R$, for every pair $(b,c) \in B_k(G_n) \times
B_k(G_n)$, there exists a path $d_1$ in $G_n$ such that\vadjust{\goodbreak} $|b| \leq R$
and $bdc$ is in $B_{2k+|d_1|}(G_n)$. We choose a single such $d_1$ for
each pair $(b,c)$, and we choose a (possibly empty) path $d_2$ such
that $bd_1cd_1$ is in $B_{2k+R(G_n)}$ (whose existence is guaranteed by
the fact that $G_n$ is irreducible). If $(b,c) \in D_{n,k}^j$, then
$bd_1cd_2$ is in $N_{n,2k+R}^{j+R}$. Using condition (C5), we have that
$R \leq m +C$ for a uniform constant $C$. Then we have that there exist
uniform constants $K_1$, $K_2$ and $K_3$ such that for each $n$, each
$k$ and each pair $(b,c)$ in $B_k(G_n) \times B_k(G_n)$,
\[
\mu\times\mu( (b,c)) \leq K_1 \lambda^{-(2m+2k)} \leq K_2 \lambda
^{-(m + R + 2k)} \leq K_3 \mu(bd_1cd_2).
\]
Thus, Lemma \ref{manyRepeatsLemma} implies that there exists a
polynomial $p_0(x)$ and $n_0$ such that for $n \geq n_0$,
\[
\mu\times\mu(D_{n,k}^j) \leq K_3 \mu(N_{n,2k+R}^{j+R}) \leq K_3
p_0(2k+R)^{\min( 2k-j, (2k+R)/U )} \lambda^{-(m+2k-j)}.
\]
With $n_1 = n_0$ and $p_1(x) = K_4 p_0(3x)^3$ for a uniform constant
$K_4$, we have
\[
\mu\times\mu(D_{n,k}^j) \leq p_1(k)^{\min( 2k-j, k/U )} \lambda^{-(m+2k-j)},
\]
which completes the proof.
\end{pf}

The following two lemmas (Lemmas \ref{DoublePerRepeatsLemma} and \ref
{manyDoublePerRepeatsLemma}) give bounds on the \mbox{$\mu\times\mu$} measure
(and therefore the cardinality) of the set of pairs of periodic paths
in $G_n$ with certain overlap properties. The general ideas are similar
to those in Lemmas \ref{manyRepeatsLemma} and \ref{masterDoubleLemma},
but in order to get precise bounds on the relevant sets, we exploit the
fact that these sets consist of pairs of periodic paths. In other
words, when we encode paths using their pattern of ``repeats,'' we also
take into account their assumed periodicity.
\begin{lemma} \label{DoublePerRepeatsLemma}
Let $(G_n)$ be a sequence of graphs satisfying the Standing Assumptions
\ref{StAssump} and bounded distortion of weights [condition
\textup{(C7)} in Definition~\ref{conditions}]. Then there exists a
polynomial $p_2(x)$ and $n_2$ such that for each $n \geq n_2$ and $k >
U(G_n)$,
\[
\mu\times\mu( S_{n,k}^{2k-1} ) \leq p_2(k) \lambda^{-(2m(G_n)+ U(G_n))}
\]
and
\[
|S_{n,k}^{2k-1}| \leq p_2(k) \lambda^{2k - U(G_n)}.
\]
\end{lemma}
\begin{pf}
Consider $(G_n)$, $n$ and $k$ as in the hypotheses. Let $m = m(G_n)$
and $U = U(G_n)$. Note that the bound on $|S_{n,k}^{2k-1}|$ follows
from the bound on $\mu\times\mu(S_{n,k}^{2k-1})$, since condition
(C7) implies that there exists a uniform constant $K$ such that
$|S_{n,k}^{2k-1}| \leq K \lambda^{2m+2k} \mu\times\mu
(S_{n,k}^{2k-1})$ [as in Lemma \ref{implicationsLemma}(1)]. We now
proceed to show the bound on $\mu\times\mu(S_{n,k}^{2k-1})$.\vspace*{1pt}

Let $b$ be in $\Per_k(G_n)$. Let $e$ be in $E_n(b)$. For $i = 1,\ldots,
k$, let $C_i \subset B_k(G_n)$ be the set of paths $c$ of length $k$ in
$G_n$ such that $c_i = e$. Then Lemma \ref{implicationsLemma} [parts
(1) and (4)] implies that there exist uniform constants $K_1$ and $K_2$
such that
%
\begin{eqnarray}\hspace*{32pt}
\mu\bigl( \Per_k(G_n) \cap C_1\bigr) &=& \mu( C_1 ) \mu(\Per_k(G_n) | C_1) \leq
K_1 \lambda^{-m} \mu(\Per_k(G_n) | C_1) \\
&\leq& K_2 \lambda^{-(m+U)}.
\end{eqnarray}
Let $C$ be the set of paths $c$ of length $k$ in $G_n$ such that $e \in
E_n(c)$. Then $C = \bigcup_{i=1}^k C_i$, and by shift-invariance of $\mu$,
%
\begin{equation} \label{aboveInequality1}
\mu\bigl( \Per_k(G_n) \cap C\bigr) \leq\sum_{i=1}^k \mu\bigl(\Per_k(G_n) \cap C_i\bigr)
\leq K_2 k \lambda^{-(m+U)}.
\end{equation}
Since $e \in E_n(b)$ was arbitrary, it follows from inequality (\ref
{aboveInequality1}) that
\begin{eqnarray*}
&&\mu\bigl( \{ c \in\Per_k(G_n)\dvtx E_n(c) \cap E_n(b) \neq\varnothing\}
\bigr)\\
&&\qquad\leq \sum_{e \in E_n(b)} \mu\bigl( \{ c \in\Per_k(G_n)\dvtx e \in E_n(c)\}\bigr) \\
&&\qquad\leq K_2 \sum_{e \in E_n(b)} k \lambda^{-(m+U)} \leq K_2 k^2 \lambda^{-(m+U)}.
\end{eqnarray*}
Since $b \in\Per_k(G_n)$ was arbitrary, we conclude that there exists
a uniform constant $K_3$ such that
\[
\mu\times\mu(S_{n,k}^{2k-1}) \leq K_2 \mu(\Per_k(G_n)) k^2 \lambda
^{-(m+U)} \leq K_3 k^2 \lambda^{-(2m+U)},
\]
where the last inequality follows from Lemma
\ref{implicationsLemma}(4). This inequality completes the proof.
\end{pf}
\begin{lemma} \label{manyDoublePerRepeatsLemma}
Let $(G_n)$ be a sequence of graphs satisfying the Standing Assumptions
\ref{StAssump} and such that $(G_n)$ has local uniqueness of paths,
small diameter and bounded distortion of weights [conditions
\textup{(C4)}, \textup{(C5)} and \textup{(C7)} in Definition
\ref{conditions}]. Then there exists a polynomial $p_3(x)$ and $n_3$
such that for $n \geq n_3$, $k > U(G_n)$ and $1 \leq j \leq2k-1$,
\[
\mu\times\mu(S_{n,k}^j) \leq p_3(k)^{k/U(G_n)} \lambda
^{-(m(G_n)+U(G_n)+2k-j)}
\]
and
\[
|S_{n,k}^j| \leq p_3(k)^{k/U(G_n)} \lambda^{j + m(G_n)-U(G_n)}.
\]
\end{lemma}
\begin{pf}
Consider $(G_n)$, $n$, $k$ and $j$ as in the hypotheses. Let
$m=m(G_n)$, $U = U(G_n)$ and $R = R(G_n)$. Note that the bound on
$|S_{n,k}^j|$ follows from the bound on $\mu\times\mu(S_{n,k}^j)$,
since condition (C7) implies that there exists a~uniform constant $K$
such that $|S_{n,k}^j| \leq K \lambda^{2m+2k} \mu\times\mu
(S_{n,k}^j)$ [as in Lem\-ma~\ref{implicationsLemma}(1)]. We now proceed
to show the bound on $\mu\times\mu(S_{n,k}^j)$.\vadjust{\goodbreak}

Let $e$ be in $E_n$ and let $C_1$ be the set of paths $b$ of length $k$
in $G_n$ such that $b_1 = e$. Then it follows from Lemma \ref
{implicationsLemma}(4) that there exists a uniform constant $K_1$ such that
%
\begin{equation} \label{ConditionalMeasureOfPeriodic}
\mu( \Per_k(G_n) | C_1) \leq K_1 \lambda^{-U}.
\end{equation}
To each\vspace*{-1pt} pair $(b,c)$ in $S_{n,k}^j$, let us associate a particular path
of length $2k+R$ in $G_n$, which we construct as follows. Let $(b,c)$
be in $S_{n,k}^j$. By definition\vspace*{1pt} of $S_{n,k}^j$, there is at least one
edge $e$ in $E_n(b) \cap E_n(c)$. Let $\tau$ be the cyclic permutation
of $\{1,\ldots, k\}$ of order $k$ given by $(1 2 \ccdots k)$. Let $\tau$
act on periodic paths of length $k$ in $G_n$ by permuting the indices:
$\tau(b_1 \ccdots b_k) = b_{\tau(1)} \ccdots b_{\tau(k)}$. Then let $b'$ be
in $\{ \tau^{\ell}(b)\dvtx \ell\in\{1,\ldots, k\}, \tau^{\ell}(b)_k = e\}
$. Similarly, let $c'$ be in $\{ \tau^{\ell}(c)\dvtx \ell\in\{1,\ldots,
k\}, \tau^{\ell}(c)_1 = e\}$. Now choose a path $d_1$ in $G_n$ such
that $|b'd_1c'| \leq R$ and $b'd_1c'$ is a path in $G_n$ (the existence
of such a path $d_1$ is guaranteed by the definition of $R$). By
irreducibility of $G_n$, we also choose a (possibly empty) path $d_2$
in $G_n$ such that $b'd_1c'd_2$ is in $B_{2k+R}$. We associate the path
$b'd_1c'd_2$ to the pair $(b,c)$, and note that there exist uniform
constants $K_2$, $K_3$ and $K_4$ [by Lemma \ref{implicationsLemma}(1)
and condition (C5)] such that
%
\begin{equation}\qquad
\mu\times\mu( (b,c) ) \leq K_2 \lambda^{-(2m +2k)}\leq K_3 \lambda
^{-(m + R + 2k)} \leq K_4 \mu(b'd_1c'd_2).
\end{equation}

Now we use the same construction as in the proof of Lemma \ref
{manyRepeatsLemma} with only slight modification. We encode the words
$b'd_1c'd_2$ as follows:
\begin{longlist}[(1)]
\item[(1)] Construct $\I$, $J$ and the partition of $\{1,\ldots, 2k+R\}$ as
in the proof of Lem\-ma~\ref{manyRepeatsLemma}, with the additional
condition that $J \cap\{k+1,\ldots, k+R\} = \varnothing$. (In other
words, we ignore any ``repeats'' introduced by $d$.)
\item[(2)] Encode $b'$ as in the proof of Lemma \ref{manyRepeatsLemma}.
\item[(3)] To encode the path $d_1$, we first encode the fact that $b'_k =
c'_1$ (by encoding $k$ and $k+|d_1|$ using an Elias code), and then
encode $d_1$ using a~prefix Shannon code with respect to $\mu( \cdot|
A_{k,k+|d_1|} \cap b')$.
\item[(4)] Encode $c'$ as in the proof of Lemma \ref{manyRepeatsLemma}.
\item[(5)] Encode $d_2$ using a prefix Shannon code with respect to $\mu(
\cdot|b'd_1c')$.
\end{longlist}
For large $n$, encoding the fact that $b'_k = c'_1$ adds less than $4
\log(2k+R)$ to $\LL(b'd_1c'd_2)$. On the other hand, we have that
there is a uniform constant $K_5 >0$ such that $\mu(A_{k,k+|d_1|} | b')
\leq K_5 \lambda^{-U}$, by Lemma \ref{implicationsLemma}(4).
Thus, there exists $n_3$ and a uniform constant $K_6$ such that for $n
\geq n_3$, we have
%
\begin{eqnarray}
&&\LL(b'd_1c'd_2) + \log\mu(b'd_1c'd_2) \nonumber\\[-8pt]\\[-8pt]
&&\qquad\leq(|\I|+1)\bigl(4 \log( 2k+R)\bigr) + N
K_6 - (m+U+|J|)\log\lambda\nonumber
\end{eqnarray}
with $|\I| \leq2k/U+2$, $N \leq|\I|$ and $|J| \geq2k-j-1$. Then by
Lemma \ref{mainInfoTheoryLemma}, there is a polynomial $p_4(x)$ such
that for $n \geq n_3$,
%
\begin{equation}
\mu\bigl( \{b'd_1c'd_2\dvtx (b,c) \in S_{n,k}^j\}\bigr) \leq p_4(k)^{k/U} \lambda
^{-(m+U + 2k - j)}.
\end{equation}
Note that the number of pairs $(b,c)$ associated to the path
$b'd_1c'd_2$ is at most $k^2$, and, hence,
%
\begin{equation}
\mu\times\mu( S_{n,k}^j) \leq k^2 p_4(k)^{k/U} \lambda^{-(m+U + 2k - j)}.
\end{equation}
Now let $p_3(x) = x^2 p_4(x)$, and the proof is complete.
\end{pf}


\subsection{Entropy}

Recall that if $G$ is a graph, then $\beta_G$ is the random variable
such that $\beta_G(\omega)$ is the spectral radius of the adjacency
matrix of $G(\omega)$.
\begin{theorem} \label{EntropyThm}
Let $(G_n)$ be a sequence of graphs that satisfies the Standing
Assumptions \ref{StAssump} and such that $(G_n)$ has local uniqueness
of paths, small diameter and bounded distortion of weights [conditions
\textup{(C4)}, \textup{(C5)} and~\textup{(C7)} in
Definition~\ref{conditions}]. Then for $1/ \lambda< \al \leq1$ and
$\varepsilon>0$,
\[
\lim_{n \rightarrow\infty} \PPP( |\beta_{G_n} - \al\lambda| \geq
\varepsilon) = 0,
\]
and the convergence to the limit is exponential in $m(G_n)$.
\end{theorem}
\begin{rmk} If we assume that $X$ is irreducible in the statement of
Theorem \ref{EntropyThmIntro}, then Theorem \ref{EntropyThmIntro} is a
direct corollary of Theorem \ref{EntropyThm}, obtained by choosing
$(G_n)$ to be the sequence of $n$-block graphs of an irreducible SFT
with positive entropy (and using the fact that such a sequence
satisfies the hypotheses of Theorem \ref{EntropyThm} by Proposition \ref
{higherGraphSeqLemma}). In the case when $X$ is reducible, $X$ has a
finite number of irreducible components of positive entropy, $X_1,\ldots,
X_r$, and there exists $i$ such that $\mathbf{h}(X_i) = \mathbf
{h}(X)$. For all large $n$, we have that $B_n(X_i) \cap B_n(X_j) =
\varnothing$ for $i \neq j$, which means that the entropies of the
random subshifts appearing inside each of these components are mutually
independent. Applying Theorem \ref{EntropyThm} to each of these
components, we obtain Theorem \ref{EntropyThmIntro} for reducible $X$.
\end{rmk}
\begin{pf*}{Proof of Theorem \ref{EntropyThm}}
Let $\al$ be in $(1/\lambda, 1]$. Let $m = m(G_n)$ and $U = U(G_n)$.
Let $b$ be a path in $G_n = (V_n,E_n)$. Let $\xi_b\dvtx \Omega_n \to\R$
be the random variable defined by
\[
\xi_b(\omega) = \cases{1, &\quad if $b$ is allowed in $G_n(\omega)$,\cr
0, &\quad else.}
\]
Now let
\[
\phi_{n,k} = \sum_{b \in B_k(G_n)} \xi_b\quad \mbox{and}\quad \psi_{n,k} =
\frac{1}{|V_n|} \sum_{b \in\Per_k(G_n)} \xi_b.
\]
For each $n$ and $k$, we have that $\psi_{n,k} \leq\beta_n^k \leq\phi
_{n,k}$. Indeed, $\psi_{n,k}$ is the average number of loops of length
$k$ based at a vertex in $G_n$. Thus, there is at least one vertex $v$
with at least $\psi_{n,k}$ loops of length $k$ based at $v$, and it
follows that $k^{-1} \log\psi_{n,k} \leq\log\beta_n$ since these
loops may be concatenated\vadjust{\goodbreak} freely. Also, it follows from subadditivity
that $\log\beta_n = \lim_k k^{-1} \log\phi_{n,k} = \inf_k k^{-1} \log
\phi_{n,k}$, which implies that $\beta_n^k \leq\phi_{n,k}$ for all $n$
and $k$.

Fix $0 < \nu< 1$, and let $k = \lceil m^{1 + \nu} \rceil+i$, where $i$
is chosen such that $0 \leq i \leq\per(G_1)-1$ and $\per(G_1)$ divides
$k$. Recall that if $(G_n)$ is the sequence of $n$-block graphs of a
fixed graph $G$, then by Proposition \ref{higherGraphSeqLemma} we have
that~$m$ and $n$ differ by at most a uniform constant, and, thus, $k
\sim n^{1+\nu}$. We will show below that as $n$ tends to infinity,
\begin{longlist}[(III)]
\item[(I)]$( \E\phi_{n,k} )^{1/k}$ tends to $\al\lambda$;
\item[(II)]$( \E\psi_{n,k} )^{1/k}$ tends to $\al\lambda$;
\item[(III)] there exists $K_1>0$ and $\rho_1>0$ such that $\frac{\Var(\phi
_{n,k})}{(\E\phi_{n,k})^2} \leq K_1 e^{-\rho_1 m}$;
\item[(IV)] there exists $K_2>0$ and $\rho_2>0$ such that $\frac{\Var(\psi
_{n,k})}{(\E\psi_{n,k})^2} \leq K_2 e^{-\rho_2 m}$.
\end{longlist}

Recall Definitions \ref{NDef}--\ref{SDef}, as well as the modification
of these definitions using ``hats.'' Notice that
\[
\E\phi_{n,k} = \sum_{b \in B_k(G_n)} \E\xi_b = \sum_{b \in B_k(G_n)}
\al^{|E_n(b)|} = \sum_{j=1}^{k} \al^j |\hat{N}^j_{n,k}|.
\]
Also,
\[
|V_n| \E\psi_{n,k} = \sum_{b \in\Per_k(G_n)} \E\xi_b = \sum_{b \in
\Per_k(G_n)} \al^{|E_n(b)|} = \sum_{j=1}^{k} \al^j |\hat{Q}^j_{n,k}|.
\]
Regarding variances, we have
\[
\Var(\phi_{n,k}) = \sum_{(b,c) \in B_k(G_n)^2} \!\al^{|E_n(b) \cup
E_n(c)|}\bigl(1 - \al^{|E_n(b) \cap E_n(c)|}\bigr)\!\leq\!\sum_{j=1}^{2k-1} \al^j
|\hat{D}_{n,k}^j|
\]
and
\[
|V_n|^2 \Var(\psi_{n,k}) = \sum_{(b,c) \in\Per_k(G_n)^2}\! \al^{|E_n(b)
\cup E_n(c)|}\bigl(1 - \al^{|E_n(b) \cap E_n(c)|}\bigr)\!\leq\!\sum_{j=1}^{2k-1} \al
^j |\hat{S}_{n,k}^j|.
\]

For the remainder of this proof, we use the following notation: if
$(x_n)$ and~$(y_n)$ are two sequences, then $x_n \sim y_n$ means that
the limit of the ratio of $x_n$ and $y_n$ tends to $1$ as $n$ tends to infinity.

\textit{Proof of} (I).
By Lemma \ref{implicationsLemma}(1), there exists a uniform constant
$K_1>0$ such that
%
\begin{equation} \label{InequalityIn(I)}\qquad\quad
\E\phi_{n,k} = \sum_{j=1}^{k} \al^j |\hat{N}^j_{n,k}| \geq\al^k \sum
_{j=1}^k |\hat{N}_{n,k}^j| = \al^k |B_k(G_n)| \geq K_1 \al^k \lambda^{m+k}.
\end{equation}
Taking $k$th roots, letting $n$ tend to infinity, and using that $m/k
\sim m^{-\nu}$ tends to $0$, we obtain that $\liminf_n ( \E\phi_{n,k}
)^{1/k} \geq\al\lambda$.\vadjust{\goodbreak}

By Lemmas \ref{implicationsLemma}(1) and \ref{manyRepeatsLemma},
we have that there exists $n_0$, a polynomial~$p_0(x)$, and a uniform
constant $K_2 >0$ such that for $n \geq n_0$,
\begin{eqnarray*}
\E\phi_{n,k} & = &\sum_{j=1}^k \al^j |\hat{N}^j_{n,k}| \\
& \leq &\sum_{j=1}^{k-1} \al^j |N^j_{n,k}| + \al^k|B_k(G_n)| \\
& \leq &( p_0(k) ) ^{k/U} \Biggl( \sum_{j=1}^{k-1} (\al\lambda)^j \Biggr) + K_2
\al^k \lambda^{k+m} \\
& \leq &(\al\lambda)^{k} \lambda^m \biggl( \frac{1}{\al\lambda-1}
p_0(k)^{k/U} \lambda^{-m} + K_2 \biggr).
\end{eqnarray*}
By condition (C4) and the fact that $k \sim m^{1+\nu}$, we have that:
\begin{itemize}
\item$m$ tends to infinity as $n$ tends to infinity by the Standing
Assumptions~\ref{StAssump};
\item$m/k \sim m^{-\nu}$, which tends to zero as $n$ tends to infinity;
\item$U \geq m -C$, which tends to infinity as $n$ tends to infinity.
\end{itemize}
Thus, taking $k$th roots and letting $n$ tend to infinity, we have that
\[
\limsup_n ( \E\phi_{n,k} )^{1/k} \leq\al\lambda,
\]
which concludes the proof of (I).

\textit{Proof of} (II).
Let $p = \per(G_1) = \per(G_n)$. Note that since $p$ divides $k$, there
exists a uniform constant $K_3 >0$ such that $|{\Per_k}(G_n)|\geq
K_3\lambda^k$ for large enough $k$. We choose $n$ large enough so that
this inequality is satisfied. Then we have that
\begin{eqnarray*}
\E\psi_{n,k} & = & |V_n|^{-1} \sum_{j=1}^{k} \al^j |\hat{Q}^j_{n,k}| \\
& \geq & |V_n|^{-1} \al^k \sum_{j=1}^k |\hat{Q}_{n,k}^j| \\
& = &|V_n|^{-1} \al^k |{\Per_k}(G_n)| \\
& \geq &K_3 \lambda^{-m}\al^k\lambda^k.
\end{eqnarray*}
Taking $k$th roots, letting $n$ tend to infinity, and using that $m/k
\sim m^{-\nu}$ tends to $0$, we get that $\liminf_n ( \E\psi_{n,k}
)^{1/k} \geq\al\lambda$. Recall that $0 \leq\psi_{n,k} \leq\phi
_{n,k}$. Therefore, it follows from (I) that $\limsup_n ( \E\psi_{n,k}
)^{1/k} \leq\al\lambda$. Thus, we have shown (II).

\textit{Proof of} (III).
For $j \leq2k-1$, Lemma \ref{masterDoubleLemma} implies that there is
$n_1$ and a~polynomial $p_1$ such that $|D_{n,k}^j| \leq p_1(k)^{k/U}
\lambda^{j+m}$ and $|D_{n,k}^{2k-1}| \leq p_1(k) \lambda^{2k+m}$ for $n
\geq n_1$. Now using that $\E\phi_{n,k} \geq K_1 \al^k \lambda^{m+k}$
[see (\ref{InequalityIn(I)})], we obtain that there exists a
uniform constant $K_5>0$ such that
\begin{eqnarray*}
\frac{\Var\phi_{n,k}}{(\E\phi_{n,k})^2} & \leq & \frac{ \sum_{j=1}^{2k-1}
\al^j |\hat{D}_{n,k}^j| }{ K_1^2 \al^{2k} \lambda^{2m+2k} } \\
& = & \frac{ \sum_{j=1}^{2k-m-1} \al^j |\hat{D}_{n,k}^j| + \sum
_{j=2k-m}^{2k-1} \al^j |\hat{D}^j_{n,k}| }{ K_1^2 \al^{2k} \lambda
^{2m+2k} } \\
& \leq & \frac{ \sum_{j=1}^{2k-m-1} \al^j |D_{n,k}^j| + \al^{2k-m} \sum
_{j=2k-m}^{2k-1} |\hat{D}^j_{n,k}| }{ K_1^2 \al^{2k} \lambda^{2m+2k} }
\\
& \leq & \frac{ p_1(k)^{k/U} \lambda^m \sum_{j=1}^{2k-m-1}(\al\lambda
)^{j} + \al^{2k-m} |D^{2k-1}_{n,k}|} { K_1^2 \al^{2k} \lambda^{2m+2k} }
\\
& \leq & \frac{ K_5 p_1(k)^{k/U} \lambda^m (\al\lambda)^{2k-m} + \al
^{2k-m} p_1(k) \lambda^{2k+m} }{ K_1^2 \al^{2k} \lambda^{2m+2k} } \\
& \leq & \frac{ K_5 }{K_1^2} \frac{p_1(k)^{k/U}}{ (\al\lambda)^m \lambda
^m } + \frac{ p_1(k) }{ K_1^2 ( \al\lambda)^m }\\
& \leq & \frac{ K_5 }{K_1^2} \biggl( \frac{p_1(k)^{k/U m}}{ (\al\lambda)}\biggr)^m
+ \frac{ p_1(k) }{ K_1^2 ( \al\lambda)^m }.
\end{eqnarray*}
Using the facts that $U \geq m - C$ and $k \sim m^{1+\nu}$, we have
that $k/Um$ is asymptotically bounded above by $2m^{\nu-1}$. Since $\nu
-1 < 0$, it holds that~$p_1(k)^{k/Um}$ tends to~$1$. Thus, we obtain
that for any $0 < \rho_1 < \ln\al\lambda$, there exists $K_6>0$ and
$n_2$ such that for $n \geq n_2$, it holds that $\Var\phi_{n,k}(\E\phi
_{n,k})^{-2} \leq K_6 e^{-\rho_1 m}$, which proves (III).

\textit{Proof of} (IV).
For $j \leq2k-1$, Lemma \ref{manyDoublePerRepeatsLemma} together with
(C4) implies that there is $n_3$ and a polynomial $p_3$ such that
$|S_{n,k}^j| \leq p_3(k)^{k/U} \lambda^j$ for $n \geq n_3$. Also,
Lem\-ma~\ref{DoublePerRepeatsLemma} implies that there is $n_4$ and a
polynomial $p_2$ such that $|S_{n,k}^{2k-1}| \leq p_2(k) \lambda
^{2k-U}$ for $n \geq n_4$. Now using that $|V_n|\E\psi_{n,k} \geq K_3
\al^k \lambda^{k}$, we obtain that there exists $K_7 >0$ such that,
with $K := K_3$,
\begin{eqnarray*}
\frac{\Var\psi_{n,k}}{(\E\psi_{n,k})^2} & \leq & \frac{ \sum_{j=1}^{2k-1}
\al^j |\hat{S}_{n,k}^j| }{ K^2 \al^{2k} \lambda^{2k} } \\
& = & \frac{ \sum_{j=1}^{2k-U-1} \al^j |\hat{S}_{n,k}^j| + \sum
_{j=2k-U}^{2k-1} \al^j |\hat{S}^j_{n,k}| }{ K^2 \al^{2k} \lambda^{2k} }
\\
& \leq & \frac{ \sum_{j=1}^{2k-U-1} \al^j |S_{n,k}^j| + \al^{2k-U} \sum
_{j=2k-U}^{2k-1} |\hat{S}^j_{n,k}| }{ K^2 \al^{2k} \lambda^{2k} } \\
& \leq & \frac{ p_3(k)^{k/U} \sum_{j=1}^{2k-U-1}(\al\lambda)^{j} + \al
^{2k-U} |S^{2k-1}_{n,k}|} { K^2 \al^{2k} \lambda^{2k} } \\
& \leq & \frac{ K_7 p_3(k)^{k/U} (\al\lambda)^{2k-U} + \al^{2k-U}
p_2(k) \lambda^{2k-U} }{ K^2 \al^{2k} \lambda^{2k} } \\
& \leq & \frac{ K_7 }{K^2} \frac{p_3(k)^{k/U}}{ (\al\lambda)^U } + \frac
{ p_2(k) }{ K^2 ( \al\lambda)^U }\\
& \leq & \frac{ K_7 }{K^2} \biggl( \frac{p_3(k)^{k/U^2}}{ (\al\lambda)}\biggr)^U +
\frac{ p_2(k) }{ K^2 ( \al\lambda)^U }.
\end{eqnarray*}
Using the facts that $U \geq m - C$ and $k \sim m^{1+\nu}$, we have
that $k/U^2$ is asymptotically bounded above by $2m^{\nu-1}$. Since $\nu
-1 < 0$, it holds that $p_3(k)^{k/U^2}$ tends to~$1$. Thus, we obtain
that for any $0 < \rho_2 < \log\al\lambda$, there exists $K_8 >0$ and
$n_5$ such that for $n \geq n_5$,
\[
\frac{\Var\phi_{n,k}}{(\E\phi_{n,k})^2} \leq K_8 e^{-\rho_2 m},
\]
which proves (IV).

\textit{Proof of Theorem} \ref{EntropyThm} \textit{using} (I)--(IV).
Recall that $\psi_{n,k} \leq\beta_n^k \leq\phi_{n,k}$. Let
$\varepsilon>0$. Since $\al\lambda>1$, we may assume without loss of
generality that \mbox{$\al\lambda- \varepsilon>1$}. Then
%
\begin{eqnarray}
\label{entropyIneq1}
&&\PPP( |\beta_n -\al\lambda| \geq\varepsilon)\nonumber\\[-8pt]\\[-8pt]
&&\qquad = \PPP( \beta_n
\geq\al\lambda+ \varepsilon) + \PPP(\beta_n \leq\al\lambda-
\varepsilon)\nonumber \\
\label{entropyIneq2}
&&\qquad = \PPP\bigl( \beta_n^k \geq(\al\lambda+ \varepsilon)^k \bigr) + \PPP\bigl( \beta
_n^k \leq(\al\lambda- \varepsilon)^k\bigr) \\
\label{entropyIneq3}
&&\qquad \leq \PPP\bigl( \phi_{n,k} \geq(\al\lambda+ \varepsilon)^k\bigr) + \PPP\bigl(\psi
_{n,k} \leq(\al\lambda- \varepsilon)^k \bigr).
\end{eqnarray}
We will bound each of the two terms in (\ref{entropyIneq3}).
Notice that
\begin{eqnarray*}
\PPP\bigl( \phi_{n,k}\!\geq\!(\al\lambda\,{+}\,\varepsilon)^k\bigr) &\,{=}\,& \PPP\bigl( \phi
_{n,k}\,{-}\,\E\phi_{n,k}\!\geq\!(\al\lambda\,{+}\,\varepsilon)^k\,{-}\,\E\phi
_{n,k}\bigr) \\
&\,{=}\,&\PPP\biggl( \phi_{n,k}\,{-}\,\E\phi_{n,k}\!\geq\!\E\phi_{n,k} \biggl(\biggl(\frac{\al
\lambda\,{+}\,\varepsilon}{(\E\phi_{n,k})^{1/k}} \biggr)^k\,{-}\,1\biggr) \biggr).
\end{eqnarray*}
Let $d^1_{n,k} = (\Var(\phi_{n,k}))^{1/2} / \E\phi_{n,k}$. Then by
Chebyshev's Inequality,
%
\begin{eqnarray}
\label{ChebBound1}
&& \PPP\bigl( \phi_{n,k}\!\geq\!(\al\lambda\,{+}\,\varepsilon)^k\bigr) \\
&&\quad\,{=}\,\PPP\biggl( \phi_{n,k}\,{-}\,\E\phi_{n,k}\!\geq\!(\Var(\phi_{n,k}))^{1/2} \frac
{1}{d^{1}_{n,k}} \biggl( \biggl(\frac{\al\lambda\,{+}\,\varepsilon}{(\E\phi
_{n,k})^{1/k}} \biggr)^k\,{-}\,1\biggr) \biggr) \hspace*{-35pt}\\
\label{ChebBound2}
&&\quad\,{\leq}\,\biggl( \frac{d^1_{n,k}}{(({\al\lambda\,{+}\,\varepsilon})/ { (\E\phi
_{n,k})^{1/k}})^k\,{-}\,1} \biggr)^2 .
\end{eqnarray}
The denominator in the right-hand side of (\ref{ChebBound2}) might be
$0$ for finitely many $n$, but by properties (I) and (III), there
exists $K_9 >0$ such that for large enough $n$,
\[
\PPP\bigl( \phi_{n,k} \geq(\al\lambda+ \varepsilon)^k\bigr) \leq\biggl( \frac
{d^1_{n,k}}{(({\al\lambda+ \varepsilon})/ { (\E\phi
_{n,k})^{1/k}})^k - 1} \biggr)^2 \leq K_9 e^{-\rho_1 m}.
\]

Similarly, we let $d^2_{n,k} = (\Var(\psi_{n,k}))^{1/2} / \E\psi
_{n,k}$, and then Chebyshev's Inequality gives that
%
\begin{eqnarray}
\label{ChebBound3}
&& \PPP\bigl( \psi_{n,k}!\leq\!(\al\lambda\,{-}\,\varepsilon)^k \bigr) \\
&&\quad\,{=}\,\PPP\biggl(
\psi_{n,k}\,{-}\,\E\psi_{n,k}\!\leq\!(\Var(\psi_{n,k}))^{1/2} \frac
{1}{d^{2}_{n,k}} \biggl( \biggl(\frac{\al\lambda\,{-}\,\varepsilon}{(\E\psi
_{n,k})^{1/k}} \biggr)^k\,{-}\,1 \biggr) \biggr) \hspace*{-35pt}\\
\label{ChebBound4}
&&\quad\,{\leq}\,\biggl( \frac{d^2_{n,k}}{(({\al\lambda\,{-}\,\varepsilon})/ { (\E\psi
_{n,k})^{1/k}})^k\,{-}\,1} \biggr)^2.
\end{eqnarray}
Again, the denominator in the right-hand side might be $0$ for finitely
many~$n$, but by properties (II) and (IV), there exists $K_{10}>0$ such
that for large enough $n$,
\[
\PPP\bigl( \psi_{n,k} \leq(\al\lambda- \varepsilon)^k \bigr) \leq\biggl( \frac
{d^2_{n,k}}{(({\al\lambda- \varepsilon})/ { (\E\psi
_{n,k})^{1/k}})^k - 1} \biggr)^2 \leq K_{10} e^{-\rho_2 m}.
\]
In conclusion, we obtain that there exists $K_{11}\,{>}\,0$ such that for
large enough~$n$,
\[
\PPP( |\beta_n -\al\lambda| \geq\varepsilon) \leq K_{11} e^{-\min(\rho
_1,\rho_2) m }.
\]
\upqed\end{pf*}

\subsection{Irreducible components of positive entropy}

\begin{theorem} \label{MSFTthm}
Let $(G_n)$ be a sequence of graphs that satisfies the Standing
Assumptions \ref{StAssump}, with $p = \per(G_1) = \per(G_n)$, and such
that:
\begin{itemize}
\item$(G_n)$ has bounded degrees [condition \textup{(C1)} in Definition \ref{conditions}],
\item$(G_n)$ has fast separation of periodic points [condition \textup{(C3)} in
Definition~\ref{conditions}],
\item and $(G_n)$ has uniform forward and backward expansion [condition
\textup{(C8)} in Definition \ref{conditions}].
\end{itemize}
Let $\mathcal{U}_{G_n}$ be the event in $\Omega_{G_n}$ that $G_n(\omega
)$ contains a unique irreducible component $C$ of positive entropy.
Also, let $\mathcal{W}_{G_n}$ be the event (contained in $\mathcal
{U}_{G_n}$) that the induced edge shift on $C$ has period $p$. Then
there exists $c>0$ such that for $1-c < \al\leq1$,
\[
\lim_{n \rightarrow\infty} \PPP( \mathcal{U}_{G_n} ) = 1\quad
\mbox{and}\quad
\lim_{n \rightarrow\infty} \PPP( \mathcal{W}_{G_n} ) = 1,
\]
and the convergence to these limits is exponential in $m(G_n)$.
\end{theorem}
\begin{rmk} Theorem \ref{MSFTthmIntro} is a corollary of Theorem \ref
{MSFTthm}: if $X$ is an irreducible SFT of positive entropy, then the
sequence of $n$-block graphs for~$X$ satisfies the hypotheses of
Theorem \ref{MSFTthm} by Proposition \ref{higherGraphSeqLemma}. In
fact, if~$X$ is a reducible SFT, we may apply Theorem \ref
{MSFTthmIntro} to each irreducible component independently, which
allows us to conclude the following. Let~$X$ be a reducible SFT with
irreducible components $X_1,\ldots, X_r$ such that $p_i = \per(X_i)$
for each $i$. Let $\mathcal{W}_n$ be the event in $\Omega_n$ that
$X_{\omega}$ has exactly $r$ irreducible components with periods
$p_1,\ldots, p_r$. Then there exists $c>0$ such that for
$\al\in(1-c,1]$, we have that $\lim_n \PPP(\mathcal{W}_n) = 1$, with
exponential (in~$n$) convergence to the limit.
\end{rmk}
%
%
\begin{defn} \label{followerDefs}
Let $G$ be a directed graph. For each vertex $v$ in $G$, and for each
$\omega$ in $\Omega_{G}$, let $\Gamma_{\omega}^+(v)$ be the union of $\{
v\}$ and the set of vertices~$u$ in $G$ such that there is an allowed
path from $v$ to $u$ in $G(\omega)$. Similarly, for each vertex $v$ in
$G$ and each $\omega$ in $\Omega_{G}$, let $\Gamma_{\omega}^-(v)$ be
the union of $\{v\}$ and the set of vertices $u$ in $G$ such that there
is an allowed path from $u$ to $v$ in $G(\omega)$. Also, let $I_{\omega
}(v) = \Gamma^+_{\omega}(v) \cap\Gamma^-_{\omega}(v)$, which is the
vertex set of the irreducible component containing $v$ in $G(\omega)$.
\end{defn}

The proof of the following proposition is an adaptation of the proof of
Lem\-ma~2.2 in \cite{ABS}.
\begin{prop} \label{posEntImpliesBig}
Let $(G_n)$ be a sequence of graphs satisfying the Standing Assumptions
\ref{StAssump} and such that $(G_n)$ has bounded degrees and uniform
forward and backward expansion [conditions \textup{(C1)} and
\textup{(C8)} in Definition \ref {conditions}]. Let $r_n$ be a sequence
of integers such that $r_n \geq a m(G_n)$, for some $a >0$, for all
large $n$. Let $C^{+}_{G_n}$ be the event in $\Omega_{G_n}$ consisting
of all $\omega$ such that there exists a vertex $v$ in $G_n$ with $r_n
\leq\Gamma_{\omega}^{+}(v) \leq |V_n|/2$. Then there exists $c > 0$
such that for $\al> 1-c$,
%
\begin{equation}
\lim_{n \rightarrow\infty} \PPP( C^{+}_{G_n} ) = 0,
\end{equation}
and the convergence of this limit is exponential in $m(G_n)$.
Furthermore, the same statement holds with \textup{``}$+$\textup{''} replaced by
\textup{``}$-$.\textup{''}
\end{prop}
\begin{pf} Let $m = m(G_n)$. Let $b >0$ be such that both $(G_n)$ and
$(G_n^T)$ are $b$-expander sequences [where the existence of such a $b$
is guaranteed by condition~(C8)]. We use the notation in Definition \ref
{followerDefs}. For any $v$ in $V_n$ and any~$\omega$ in $\Omega
_{G_n}$, the set $\Gamma_{\omega}^+(v)$ has the property that all edges
in $E_n(\Gamma_{\omega}^+(v),\overline{\Gamma_{\omega}^+(v)})$ are
forbidden (by $\omega$). Then the fact that $G_n$ is a $b$-expander
implies that for a particular subset $S$ of $V_n$, the probability that
$S = \Gamma_{\omega}^+(v)$ for some $v$ is bounded above by $(1-\al)^{b
|S|}$. The number of subsets $S$ of $V_n$ with $|S|=r$ that could
appear as $\Gamma_{\omega}^+(v)$ for some $v$ is bounded above by
$(\Delta e)^r$, where $e$ is the base of the natural\vadjust{\goodbreak} logarithm~(\cite
{ABS}, Lemma 2.2) (see also \cite{Alon}, Lemma~2.1, or \cite{Knuth},
page 396, Exercise~11). Then for $\al$ such that $\Delta e (1-\al) <
1$, we have that for any $0 \leq r_n \leq|V_n|/2$,
%
\begin{eqnarray}
\PPP(C^{+}_{G_n}) & = & \PPP\biggl( \exists v \mbox{ such that } r_n \leq
|\Gamma_{\omega}^+(v)| \leq\frac{|V_n|}{2}\biggr) \\
&\leq& \sum
_{r=r_n}^{{|V_n|}/{2}} |V_n| (\Delta e)^r (1-\al)^{b r} \\
&\leq& |V_n|\bigl(\Delta e (1-\al)^b\bigr)^{r_n} \frac{1}{1-\Delta e (1-\al)} \\
\label{RHSgoesToZero}
&\leq& \bigl(\lambda^{1/a} \Delta e (1-\al)^b\bigr)^{a m} \frac{1}{1-\Delta e
(1-\al)}.
\end{eqnarray}
Thus, there is a $c>0$ (depending only on $a$, $b$, $\lambda$ and
$\Delta$) such that if $\al> 1-c$, then the right-hand side of the
inequality in (\ref{RHSgoesToZero}) tends to zero exponentially in
$m(G_n)$ as $n$ tends to infinity. In particular, we may take
\[
c = \biggl( \frac{1}{\lambda} \biggr)^{1/{ab}} \biggl( \frac{1}{\Delta e} \biggr)^{1/b}.
\]

Since $(G_n^T)$ is also a uniform $b$-expander, the same estimates hold
with~$C^{-}_{G_n}$ in place of $C^{+}_{G_n}$.
\end{pf}
\begin{pf*}{Proof of Theorem \ref{MSFTthm}}
Let $(G_n)$ be as in the statement of Theorem~\ref{MSFTthm}. Let $m =
m(G_n)$, $z = z(G_n)$, and $p = \per(G_1)=\per(G_n)$. We use the
notation in Definition \ref{followerDefs}. Consider the following events:
\begin{eqnarray*}
F^+_n & = & \{ \omega\in\Omega_n\dvtx \forall v \in V_n, \Gamma^+_{\omega
}(v) \leq z(G_n) -2p \mbox{ or } \Gamma^+_{\omega}(v) > |V_n|/2 \} ,\\
F^-_n & = & \{ \omega\in\Omega_n\dvtx \forall v \in V_n, \Gamma^-_{\omega
}(v) \leq z(G_n) -2p \mbox{ or } \Gamma^-_{\omega}(v) > |V_n|/2 \}, \\
F_n & = & F_n^+ \cap F_n^-.
\end{eqnarray*}
Recall that condition (C3) gives $a>0$ such that $z \geq a m$. Note
that Proposition \ref{posEntImpliesBig} gives $c > 0$ such that for
$1-c < \al\leq1$, there exists $K_1, K_2>0$ and $\rho_1, \rho_2 >0$
such that for large $n$,
\[
\PPP(\Omega_n \setminus F_n^+) \leq K_1 e^{-\rho_1 m} \quad\mbox{and}\quad \PPP
(\Omega_n \setminus F^-_n) \leq K_2 e^{-\rho_2 m}.
\]
Fix such an $\al$, and note that for all large enough $n$, we have the
following estimate: $\PPP(\Omega_n \setminus F_n) \leq2\max(K_1,K_2)
e^{-\min(\rho_1,\rho_2)m}$.

Consider $\omega$ in $F_n$. Suppose that there exists $v_1$ and $v_2$
in $V_n$ such that $|I_{\omega}(v_1)| > z -2p$ and $|I_{\omega}(v_2)| >
z -2p$. Then by definition of $F_n$, we must have that $\Gamma^+_{\omega
}(v_1) \cap\Gamma^-_{\omega}(v_2) \neq\varnothing$ and $\Gamma
^-_{\omega}(v_1) \cap\Gamma^+_{\omega}(v_2) \neq\varnothing$. It
follows that there is a path from $v_1$ to $v_2$ in $G_n(\omega)$, and
there is a path from $v_2$ to $v_1$ in $G_n(\omega)$. Thus, $I_{\omega
}(v_1) = I_{\omega}(v_2)$. We have shown that for $\omega$ in $F_n$,
there is at most one irreducible component of cardinality greater than
$z-2p$.\vadjust{\goodbreak} Note that this argument implies that for $\omega$ in $F_n$, all
allowed periodic orbits $\gamma$ such that $|V_n(\gamma)|>z-2p$ must
lie in the same irreducible component.

By definition of $z$, if $I_{\omega}$ is an irreducible component of
$G_n(\omega)$ with positive entropy, then $|I_{\omega}| > z$ (since it
must contain at least two periodic orbits with overlapping vertex
sets). We deduce that for $\omega$ in $F_n$, there is at most one
irreducible component of $G_n(\omega)$ with positive entropy.

We now show that there exists an irreducible component of positive
entropy with probability tending exponentially to $1$. Let $z_1 = z-i$,
where $i$ is chosen (for each $n$) such that $0 \leq i \leq p-1$ and
$p$ divides $z_1$. Then let $z_2 = z_1-p$. Consider the following
sequences of random variables:
%
\begin{equation}
f_n = \sum_{b \in\Per_{z_1}(G_n)} \xi_b\quad \mbox{and}\quad g_n = \sum_{b
\in\Per_{z_2}(G_n)} \xi_b.
\end{equation}
Note that by the definition of $z$ and Lemma \ref{perLemma}, we have
that $|E_n(b)| = |b|$ for any periodic path $b$ with period less than
or equal to $z$. Furthermore, any two such paths are disjoint.
Therefore, the random variables $\{ \xi_b\}_{b \in\Per_{z_1}(G_n)}$
are jointly independent, and the random variables $\{ \xi_b\}_{b \in
\Per_{z_2}(G_n)}$ are also jointly independent. Thus,
\begin{eqnarray*}
\E f_n & = & \sum_{b \in\Per_z(G_n)} \al^{z_1} = \al^{z_1} |{\Per
_{z_1}}(G_n)|, \\
\E g_n & = & \sum_{b \in\Per_{z_2}(G_n)} \al^{z_2} = \al^{z_2}|{\Per
_{z_2}}(G_n)|, \\
\Var(f_n) & = & \sum_{b \in\Per_{z_1}(G_n)} \al^{z_1}(1-\al^{z_1}) = \al
^{z_1} (1-\al^{z_1}) |{\Per_{z_1}}(G_n)|, \\
\Var(g_n) & = & \sum_{b \in\Per_{z_2}(G_n)} \al^{z_2}(1-\al^{z_2}) = \al
^{z_2} (1-\al^{z_2}) |{\Per_{z_2}}(G_n)|.
\end{eqnarray*}
As $n$ tends to infinity, $z$ tends to infinity since $z \geq a m$ and
$m$ tends to infinity. Then by the Standing Assumptions \ref{StAssump}
[in particular, we use that $\Sp_{\times}(G_n) = \Sp_{\times}(G_1)$]
and the fact that $p$ divides $z_1$ and $z_2$, we have that each of the
sequences $\lambda^{-z_1}|{\Per_{z_1}}(G_n)|$ and $\lambda^{-z_2} |{\Per
_{z_2}}(G_n)|$ tends to a finite, nonzero limit as $n$ tends to infinity
(and in fact the limit is $p$). For two sequences $x_n$ and $y_n$ of
positive real numbers, let $x_n \sim y_n$ denote the statement that
their ratio tends to a finite, nonzero limit as $n$ tends to infinity.
Then we have that $\E f_n \sim(\al\lambda)^{z_1} \sim\Var(f_n)$ and
$\E g_n \sim(\al\lambda)^{z_2} \sim\Var(g_n)$. Note that $\E f_n
\geq\Var(f_n)$ and $\E g_n \geq\Var(g_n)$. A~simple application of
Chebyshev's Inequality implies that
\begin{eqnarray*}
\PPP( f_n \leq0 ) & \leq& \PPP\bigl( f_n - \E f_n \leq- \Var(f_n) \bigr) \\
& \leq & \biggl( \frac{1}{\Var(f_n)^{1/2}} \biggr)^2 \sim\biggl( \frac{1}{\al\lambda}
\biggr)^{z_1} \leq\biggl( \frac{1}{\al\lambda} \biggr)^{a m - i}
\end{eqnarray*}
and
\begin{eqnarray*}
\PPP( g_n \leq0 ) & \leq & \PPP\bigl( g_n - \E g_n \leq- \Var(g_n) \bigr) \\
& \leq & \biggl( \frac{1}{\Var(g_n)^{1/2}} \biggr)^2 \sim\biggl( \frac{1}{\al\lambda}
\biggr)^{z_2} \leq\biggl( \frac{1}{\al\lambda} \biggr)^{a m-i-p}.
\end{eqnarray*}
We have shown that the probability that there is no periodic orbit of
period~$z_1$ tends to $0$ exponentially in $m$ as $n$ tends to
infinity, and the probability that there exists no periodic orbit of
period $z_2$ tends to $0$ exponentially in~$m$ as $n$ tends to infinity.

In summary, we have shown that the following events occur with
probability tending to $1$ exponentially in $m$ as $n$ tends to infinity:
\begin{itemize}
\item there exists a periodic point of period $z-i$;
\item there exists a periodic point of period $z-i-p$;
\item any two periodic points of period greater than $z - 2p$ lie in
the same irreducible component (of necessarily positive entropy);
\item there is at most one irreducible component of positive entropy.
\end{itemize}
We conclude that with probability tending to $1$ exponentially in $m$
as $n$ tends to infinity, there exists a unique irreducible component
of positive entropy, and the induced edge shift on that component has
period $p$.
\end{pf*}

\section{Remarks} \label{Remarks}

\begin{rmk} \label{relaxRmk}
The proofs of Theorems \ref{EmptyThm} and \ref{SubCritThm} do not
require all of the Standing Assumptions \ref{StAssump}. In fact, these
proofs only use that $\Sp_{\times}(G_n) = \Sp_{\times}(G_1)$ for each
$n$ and that $z(G_n)$ tends to infinity as $n$ tends to infinity.
\end{rmk}
\begin{rmk}
Theorem \ref{EmptyThm} states that at the critical threshold $\al=
1/\lambda$, the probability of emptiness tends to zero. Using the fact
that entropy is a~monotone increasing random variable (as defined in
Section \ref{ProbFramework}), one may deduce from Theorem \ref
{EntropyThm} that for $\al= 1/\lambda$, the probability that the
random SFT has zero entropy tends to $1$. It might be interesting to
know more about the behavior of typical random SFTs at the critical threshold.
\end{rmk}
\begin{rmk}
We have considered only random $\Z$-SFTs, but one may also consider
random $\Z^d$-SFTs for any $d$ in $\N$ by adapting the construction of~%
$\Omega_n$ and $\PPP$ in the obvious way. It appears that most of the
proofs presented above may not be immediately adapted for $d>1$, but
there is one exception, which we state below. Let $X$ be a nonempty $\Z
^d$-SFT. For $d >1$, there are various zeta functions for $X$ (for a
definition distinct from ours, see \cite{Lind1}); we consider
\[
\zeta_{X}(t) = \exp\Biggl( \sum_{p=1}^{\infty} \frac{N_p}{p} t^p \Biggr),
\]
where $N_p$ is the number of periodic points $x$ in $X$ such that the
number of points in the orbit of $x$ divides $p$. The function $\zeta
_X$ has radius of convergence $1/\rho$, where $\log\rho= \limsup_p
p^{-1} \log(N_p)$. For example, for a full $\Z^d$ shift on $a$ symbols,
\mbox{$\rho= a$}, regardless of $d$. Using exactly the same proof as
presented in Section \ref{Emptiness}, we obtain that
\[
\limsup_{n \to\infty} \PPP(\mathcal{E}_n) \leq\cases{(\zeta_{X}(\al
))^{-1}, &\quad if $\al\in[0, 1/\rho)$,\cr
0, &\quad if $\al\in[1/\rho,1]$.}
\]
For $\al\geq1/\rho$, this bound implies that the limiting probability
of emptiness is $0$. In this context, we note that there is no
algorithm, which, given a $\Z^d$-SFT $X$ defined by a finite list of
finite forbidden configurations, will decide whether $X$ is empty~\cite
{Berger}. Nonetheless, we may be able to compute the limiting
probability of emptiness. For example, if $X$ is a full shift on $a$
symbols, then for $\al\geq1/a$, we have that the limiting probability
of emptiness is $0$.\vspace*{-2pt}
\end{rmk}
\begin{rmk}
One may also consider more general random subshifts. Recall that a set
$X \subset\cA^{\Z}$ is a subshift if it is closed and shift-invariant.
For a nonempty subshift $X$ and a natural number $n$, we may consider
the (finite) set of subshifts obtained by forbidding words of length
$n$ from $X$. After defining a probability measure $\PPP$ on this space
as in Section \ref{Preliminaries}, we obtain random subshifts of $X$.
Now we may investigate the asymptotic probability of properties of
these random subshifts. Recall that any subshift $X$ can be written as
$\bigcap X_n$, where $(X_n)$ is a sequence of SFTs (called the Markov
approximations of $X$) and $\lim_n h(X_n) = h(X)$. A subshift $X$ is
called \textit{almost sofic} \cite{Petersen1} if there exists a
sequence $(X_n)$ of irreducible SFTs such that \mbox{$X_n \subset X$} and $\lim
_n h(X_n) = h(X)$. Using this inner and outer approximation by SFTs,
the conclusion of Theorem \ref{EntropyThmIntro} still holds if the
system $X$ is only assumed to be an almost sofic subshift.\vspace*{-2pt}
\end{rmk}

\begin{rmk}
Theorem \ref{MSFTthm} asserts the existence of a constant \mbox{$c>0$}, but we
are left with several questions about this constant. Fix a sequence
$(G_n)$ satisfying the hypotheses of Theorem \ref{MSFTthm}. Let $\al_*
= \inf\{\al>0\dvtx \lim_n \PPP(\mathcal{U}_n) =1\}$. What is~$\al_*$?
What is $\al_*$ in the case that $(G_n)$ is the sequence of $n$-block
graphs of a mixing SFT of positive entropy (or even a~full shift)?\vspace*{-2pt}
\end{rmk}


\section*{Acknowledgments}
The author expresses great thanks to his advisor, Mike Boyle, for his
guidance and encouragement in general, and also for many helpful
conversations and comments regarding this work in particular. Further,
the author thanks the anonymous referee for reading the paper carefully
and making helpful suggestions.\vspace*{-2pt}



%
\printaddresses

\end{document}